%% file: arXiv_revised.tex
\documentclass[letterpaper,11pt]{article}
\usepackage[margin=1in]{geometry}  
\input{defs}
\input{preamble.tex}
\begin{document}
\title{Sharp mean-field analysis of permutation mixtures\\
and permutation-invariant decisions}
\author{Yiguo Liang, Yanjun Han\thanks{Yiguo Liang is with New York University, Shanghai, China. Yanjun Han is with the Courant Institute of Mathematical Sciences and Center for Data Science, New York University, New York, USA. Emails: \url{{yl10351, yanjunhan}@nyu.edu}. }}
\maketitle

\begin{abstract} 
We develop sharp bounds on the statistical distance between high-dimensional permutation mixtures and their i.i.d. counterparts. Our approach establishes a new geometric link between the spectrum of a complex channel overlap matrix and the information geometry of the channel, yielding tight dimension-independent bounds that close gaps left by previous work. Within this geometric framework, we also derive dimension-dependent bounds that uncover phase transitions in dimensionality for Gaussian and Poisson families. Applied to compound decision problems, this refined control of permutation mixtures enables sharper mean-field analyses of permutation-invariant decision rules, yielding strong non-asymptotic equivalence results between two notions of compound regret in Gaussian and Poisson models.
\end{abstract}

\tableofcontents

\input{Intro.tex}

\input{Improved_Bounds_for_Permutation_Mixtures.tex}

\input{Application_To_Compound_Decision_Problems.tex}

\input{Discussion.tex}
 
\appendix
\input{Deferred_Proofs.tex }

\bibliographystyle{alpha}
\bibliography{ref.bib}

\end{document}

%% file: defs.tex
\usepackage{bbm}
\usepackage{graphicx}
\usepackage{amsmath,amssymb,amsthm,amsfonts}

\usepackage{paralist}
\usepackage{bm}
\usepackage{xspace}
\usepackage{url}
\usepackage{prettyref}
\usepackage{boxedminipage}
\usepackage{wrapfig}
\usepackage{ifthen}
\usepackage{color}
\usepackage{xspace}

\usepackage{amsmath,amsthm,amsfonts,amssymb}
\usepackage{graphicx}

\usepackage{nicefrac}

\newtheorem*{definition*}{Definition}

\usepackage{subcaption}

\usepackage[utf8]{inputenc}

\usepackage{xcolor}
\definecolor{expert}{HTML}{008000}
\definecolor{error}{HTML}{f96565}

\usepackage{color-edits}
\addauthor{sw}{blue}

\usepackage{thmtools}
\usepackage{thm-restate}

\usepackage{tikz}
\usetikzlibrary{arrows,calc} 
\newcommand{\tikzAngleOfLine}{\tikz@AngleOfLine}
\def\tikz@AngleOfLine(#1)(#2)#3{%
\pgfmathanglebetweenpoints{%
\pgfpointanchor{#1}{center}}{%
\pgfpointanchor{#2}{center}}
\pgfmathsetmacro{#3}{\pgfmathresult}%
}

\declaretheoremstyle[
    headfont=\normalfont\bfseries, 
    bodyfont = \normalfont\itshape]{mystyle} 

\usepackage{listings}
\usepackage{amsmath}
\usepackage{amsthm}
\usepackage{tikz}
\usepackage{caption}
\usepackage{mdwmath}
\usepackage{multirow}
\usepackage{mdwtab}
\usepackage{eqparbox}
\usepackage{multicol}
\usepackage{amsfonts}
\usepackage{tikz}
\usepackage{multirow,bigstrut,threeparttable}
\usepackage{amsthm}
\usepackage{bbm}
\usepackage{epstopdf}
\usepackage{mdwmath}
\usepackage{mdwtab}
\usepackage{eqparbox}
\usetikzlibrary{topaths,calc}
\usepackage{latexsym}
\usepackage{cite}
\usepackage{amssymb}
\usepackage{bm}
\usepackage{amssymb}
\usepackage{graphicx}
\usepackage{mathrsfs}
\usepackage{epsfig}
\usepackage{psfrag}
\usepackage{setspace}
\usepackage[
            CJKbookmarks=true,
            bookmarksnumbered=true,
            bookmarksopen=true,
            colorlinks=true,
            citecolor=red,
            linkcolor=blue,
            anchorcolor=red,
            urlcolor=blue
            ]{hyperref}

\usepackage{comment}
\usepackage{mathtools}
\usepackage{blkarray}
\usepackage{multirow,bigdelim,dcolumn,booktabs}

\usepackage{xparse}
\usepackage{tikz}
\usetikzlibrary{calc}
\usetikzlibrary{decorations.pathreplacing,matrix,positioning}

\usepackage[T1]{fontenc}
\usepackage[utf8]{inputenc}
\usepackage{mathtools}
\usepackage{blkarray, bigstrut}
\usepackage{gauss}

\newcommand*{\BraceAmplitude}{0.4em}%
\newcommand*{\VerticalOffset}{0.5ex}%
\newcommand*{\HorizontalOffset}{0.0em}%
\newcommand*{\blocktextwid}{3.0cm}%
\NewDocumentCommand{\InsertLeftBrace}{%
	O{} 
	O{\HorizontalOffset,\VerticalOffset} 
	O{\blocktextwid} 
	m   
	m   
	m   
}{%
	\begin{tikzpicture}[overlay,remember picture]
	\coordinate (Brace Top)    at ($(#4.north) + (#2)$);
	\coordinate (Brace Bottom) at ($(#5.south) + (#2)$);
	\draw [decoration={brace, amplitude=\BraceAmplitude}, decorate, thick, draw=black, #1]
	(Brace Bottom) -- (Brace Top) 
	node [pos=0.5, anchor=east, align=left, text width=#3, color=black, xshift=\BraceAmplitude] {#6};
	\end{tikzpicture}%
}%
\NewDocumentCommand{\InsertRightBrace}{%
	O{} 
	O{\HorizontalOffset,\VerticalOffset} 
	O{\blocktextwid} 
	m   
	m   
	m   
}{%
	\begin{tikzpicture}[overlay,remember picture]
	\coordinate (Brace Top)    at ($(#4.north) + (#2)$);
	\coordinate (Brace Bottom) at ($(#5.south) + (#2)$);
	\draw [decoration={brace, amplitude=\BraceAmplitude}, decorate, thick, draw=black, #1]
	(Brace Top) -- (Brace Bottom) 
	node [pos=0.5, anchor=west, align=left, text width=#3, color=black, xshift=\BraceAmplitude] {#6};
	\end{tikzpicture}%
}%
\NewDocumentCommand{\InsertTopBrace}{%
	O{} 
	O{\HorizontalOffset,\VerticalOffset} 
	O{\blocktextwid} 
	m   
	m   
	m   
}{%
	\begin{tikzpicture}[overlay,remember picture]
	\coordinate (Brace Top)    at ($(#4.west) + (#2)$);
	\coordinate (Brace Bottom) at ($(#5.east) + (#2)$);
	\draw [decoration={brace, amplitude=\BraceAmplitude}, decorate, thick, draw=black, #1]
	(Brace Top) -- (Brace Bottom) 
	node [pos=0.5, anchor=south, align=left, text width=#3, color=black, xshift=\BraceAmplitude] {#6};
	\end{tikzpicture}%
}%

\usetikzlibrary{patterns}

\definecolor{cof}{RGB}{219,144,71}
\definecolor{pur}{RGB}{186,146,162}
\definecolor{greeo}{RGB}{91,173,69}
\definecolor{greet}{RGB}{52,111,72}

\theoremstyle{plain}
\newtheorem{thm}{Theorem}[section]

\newtheorem{lemma}[thm]{Lemma}
\newtheorem{remark}[thm]{Remark}
\newtheorem{corollary}[thm]{Corollary}
\newtheorem{definition}[thm]{Definition}

\def \bP {\mathbb{P}}
\def \bQ {\mathbb{Q}}

\def \bE {\mathbb{E}}
\def \bR {\mathbb{R}}

\def \var {\mathrm{Var}}
\def\1{\mathbbm{1}}

\usepackage{xspace}

\newcommand{\stepa}[1]{\overset{\rm (a)}{#1}}
\newcommand{\stepb}[1]{\overset{\rm (b)}{#1}}
\newcommand{\stepc}[1]{\overset{\rm (c)}{#1}}
\newcommand{\stepd}[1]{\overset{\rm (d)}{#1}}

\newcommand{\naturals}{\mathbb{N}}

\newcommand{\TV}{{\sf TV}}

\newcommand{\KL}{{\sf KL}}

\newcommand{\pth}[1]{\left( #1 \right)}
\newcommand{\qth}[1]{\left[ #1 \right]}
\newcommand{\sth}[1]{\left\{ #1 \right\}}
\newcommand{\bpth}[1]{\Big( #1 \Big)}
\newcommand{\bqth}[1]{\Big[ #1 \Big]}
\newcommand{\bsth}[1]{\Big\{ #1 \Big\}}

\newcommand{\Poi}{\text{\rm Poi}}
\newcommand{\Unif}{\text{\rm Unif}}

\newcommand{\indc}[1]{{\mathbf{1}_{\left\{{#1}\right\}}}}

\definecolor{myblue}{rgb}{.8, .8, 1}
\definecolor{mathblue}{rgb}{0.2472, 0.24, 0.6} 
\definecolor{mathred}{rgb}{0.6, 0.24, 0.442893}
\definecolor{mathyellow}{rgb}{0.6, 0.547014, 0.24}

\newcommand{\blue}{}

\newcommand{\sfC}{{\mathsf{C}}}
\newcommand{\sfD}{{\mathsf{D}}}

\newcommand{\calC}{{\mathcal{C}}}
\newcommand{\calD}{{\mathcal{D}}}

\newcommand{\calN}{{\mathcal{N}}}

\newcommand{\calP}{{\mathcal{P}}}

\newcommand{\calT}{{\mathcal{T}}}

\newcommand{\rmd}{\mathrm{d}}
\newcommand{\Perm}{\mathrm{Perm}}

\newcommand{\trace}{\mathrm{Tr}}

\usepackage{cleveref}
\crefname{lemma}{Lemma}{Lemmas}
\Crefname{lemma}{Lemma}{Lemmas}
\crefname{thm}{Theorem}{Theorems}
\Crefname{thm}{Theorem}{Theorems}
\Crefname{assumption}{Assumption}{Assumptions}
\crefformat{equation}{(#2#1#3)}

\newcommand{\Capa}{\sfC_{\chi^2}(\calP)}
\newcommand{\Diam}{\sfD_{k}(\calP)}

\newcommand{\MSE}{{\sf MSE}}

\newcommand{\thetahatS}{{\widehat{\theta}^{\mathrm{S}}}}
\newcommand{\thetahatPI}{{\widehat{\theta}^{\mathrm{PI}}}}
\newcommand{\regS}{\mathrm{TotReg}^{\mathrm{S}}}
\newcommand{\regPI}{\mathrm{TotReg}^{\mathrm{PI}}}

%% file: preamble.tex
\usepackage{stfloats}
\usepackage{cuted}
\allowdisplaybreaks

%% file: Intro.tex
\section{Introduction}
  Mean-field approximation is a powerful framework for studying the behavior of large random systems with weak dependencies. It enables precise asymptotic computations of key quantities such as the free energy in complex models arising from probability theory and statistical physics \cite{Par88,JaiRisKoe19,LacMukYeu24}. A central theoretical question is to quantify when a high-dimensional distribution can be accurately approximated by a product measure, i.e., to rigorously understand the effectiveness of mean-field approximations. In this work, we study mean-field approximations for two connected high-dimensional problems: \emph{permutation mixtures} in probability theory and \emph{compound decision problems} in statistics. We begin by recalling the notion of permutation mixtures, introduced in the recent work \cite{han2024approximate}. 
\begin{definition}[Permutation mixture and its mean-field approximation]\label{def:perm_mixture}
Let $\calP$ be a collection of probability measures on a fixed probability space.
Given $P_1,\cdots,P_n\in \calP$: 
\begin{itemize}
    \item The distribution $\bP_n$ of the observation $X=(X_1,\cdots,X_n)$ is the ``permutation mixture'' defined by
    \begin{align}\label{eq:bP_n}
        (X_1,\cdots,X_n) \sim \bE_{\pi\sim \Unif(S_n)}\qth{ \otimes_{i=1}^n P_{\pi(i)} }; 
    \end{align}
    \item The distribution $\bQ_n$ of the observation $X=(X_1,\cdots,X_n)$ is an i.i.d.\ product of one-dimensional mixtures defined by
    \begin{align}\label{eq:bQ_n}
        (X_1,\cdots,X_n) \sim \pth{\frac{1}{n}\sum_{i=1}^n P_i}^{\otimes n}. 
    \end{align}
\end{itemize}
\end{definition}
\par\noindent
 In words, the permutation mixture $\bP_n$ is the joint distribution of a uniformly random permutation of a random vector with independent coordinates, and the mean-field approximation $\bQ_n$ is a product distribution with the same one-dimensional marginals as $\bP_n$. A key discovery of \cite{han2024approximate} (see also the earlier work \cite{Din22,tang2023capacity} in special cases detailed in \Cref{subsec:related-work}) is that under mild conditions on $\calP$ (where no pair of distributions in $\calP$ is mutually singular), the statistical distance between $\bP_n$ and $\bQ_n$, albeit being two high-dimensional distributions, remains bounded \emph{regardless of the dimension $n$}. As a specific example, when $\calP = \{ \calN(\theta,1): |\theta|\le \mu \}$ is the normal mean model, \cite[Corollary 1.3]{han2024approximate} shows that
\begin{align}\label{eq:old_upper_bound}
\sup_n \sup_{P_1,\dots,P_n\in \calP} \chi^2\pth{ \bP_n \| \bQ_n } = \begin{cases}
    O(\mu^4) &\text{if } \mu \le 1, \\
    \exp(O(\mu^3)) & \text{if } \mu > 1. 
\end{cases}
\end{align}
Using the inequality $2\TV(\bP_n, \bQ_n)^2\le \KL(\bP_n \| \bQ_n) \le \log(1+\chi^2(\bP_n \| \bQ_n))$, this yields upper bounds on the total variation (TV) distance and Kullback--Leibler (KL) divergence as well. Lower bounds of similar form were also established in \cite[Lemma 6.1]{han2024approximate}: 
\begin{align}\label{eq:old_lower_bound}
\sup_n \sup_{P_1,\dots,P_n\in \calP} \chi^2\pth{ \bP_n \| \bQ_n } = \begin{cases}
    \Omega(\mu^4) &\text{if } \mu \le 1, \\
    \exp(\Omega(\mu^2)) & \text{if } \mu > 1, 
\end{cases}
\end{align}
demonstrating that the phase transition of the $\chi^2$ divergence as $\mu$ crosses $1$ is real. 

However, the upper bound in \eqref{eq:old_upper_bound} still leaves significant gaps when $\mu>1$. First, \cite{han2024approximate} conjectured that the upper bound could be improved to $\exp(O(\mu^2))$, matching the lower bound \eqref{eq:old_lower_bound}. Achieving this improvement is technically challenging because the analysis in \cite{han2024approximate} depends on all eigenvalues of a complicated channel overlap matrix (see \Cref{defn:channel-overlap-matrix}); while \cite{han2024approximate} sidestepped this complexity by using only the trace and spectral gap, such relaxations may be fundamentally suboptimal. Second, even if the improved bound $\exp(O(\mu^2))$ holds, it indicates an incorrect dependence on $\mu$ unless the dimension $n$ is \emph{extremely large}. In fact, we will show that (see \Cref{cor:dim-dep-bound}) it only becomes tight when $n > \exp(\Omega(\mu^2))$, an extremely large threshold for large $\mu$; below this threshold, the correct scaling of the $\chi^2$ divergence is $\exp(\Theta(\mu\sqrt{\log n}))$, exhibiting a linear dependence on $\mu$ in the exponent. These observations motivate the need for new techniques to provide tight mean-field approximations for general permutation mixtures, in both \emph{dimension-independent} and \emph{dimension-dependent} scenarios. 

While these improvements may initially appear technical, they are in fact essential for making progress on compound decision problems, a classical topic introduced by Robbins \cite{robbins1951asymptotically}. A compound decision problem consists of $n$ independent estimation tasks with shared structure. Formally, let $\{P_\theta: \theta\in \Theta \subseteq \bR\}$ be a parametric family of distributions, the learner observes $X_i\sim P_{\theta_i}$ independently for $i\in [n]$ and is tasked to estimate the unknown parameters $\theta_1,\dots,\theta_n$. An estimator (or decision rule) $\widehat{\theta}$ is a mapping from observations $(X_1,\dots,X_n)$ to estimates $(\widehat{\theta}_1,\dots,\widehat{\theta}_n)$, with mean squared error (MSE) given by
\begin{align}\label{eq:MSE}
\MSE(\theta, \widehat{\theta}) = \bE_{\theta} \qth{\| \widehat{\theta}(X^n) - \theta \|^2}. 
\end{align}
To solve this estimation problem, the empirical Bayes (EB) framework \cite{robbins1951asymptotically,robbins1956empirical} \emph{pretends} that $\theta_1,\dots,\theta_n$ were i.i.d. drawn from the (unknown) empirical distribution $G_n = \frac{1}{n}\sum_{i=1}^n \delta_{\theta_i}$, and competes with the Bayes estimator under this prior, a.k.a. the \emph{separable/simple oracle}. Equivalently, the separable oracle $\thetahatS$ is the optimal estimator in the class of coordinate-wise (separable) decision rules (denote by $\calD^{\rm S}$): $\widehat{\theta}(X_1,\dots,X_n) = (f(X_1),\dots,f(X_n))$ for some function $f$, and satisfies
\begin{align}\label{eq:separable_oracle}
\MSE(\theta, \thetahatS) \le \MSE(\theta, \widehat{\theta}), \quad \forall \widehat{\theta}\in \calD^{\rm S}. 
\end{align}
Even though $G_n$ is unknown, EB procedures propose to mimic the Bayes estimator, either explicitly (by estimating the prior and applying the learned Bayes estimator) or implicitly (by approximating the Bayes estimator directly). The appealing aspect of Robbins's EB framework is that, the resulting EB estimator $\widehat{\theta}^{\rm EB}$ often satisfies
\begin{align}\label{eq:regret_against_simple}
\regS(\widehat{\theta}^{\rm EB}) &:= \sup_{\theta\in \Theta} \pth{ \MSE(\theta,\widehat{\theta}^{\rm EB}) - \MSE(\theta,\thetahatS) }= o(n),
\end{align}
showing that coordinatewise, $\widehat{\theta}^{\rm EB}$ asymptotically achieves the same MSE at
the unknown instance $\theta$ as the best separable rule that knows $\theta$. The quantity $\regS(\widehat{\theta})$ is called the \emph{total regret} of $\widehat{\theta}$ against the separable oracle, and characterizing its rate of convergence has been the central focus of a series of papers \cite{robbins1951asymptotically,Jiang_2009,brown2009nonparametric,polyanskiy2021sharp,shen2022empirical,jana2022optimal,jana2023empirical,ghosh2025stein}.

However, the EB estimator is not typically separable, so \eqref{eq:separable_oracle} does not hold for $\widehat{\theta} = \widehat{\theta}^{\rm EB}$ and the separable oracle is not always a compelling benchmark. By contrast, the EB estimator typically belongs to a larger class of estimators (denoted by $\calD^{\rm PI}$): the \emph{permutation-invariant} rules, satisfying
\begin{align*}
\widehat{\theta}_{\pi(i)}(X_{\pi(1)},\dots,X_{\pi(n)}) = \widehat{\theta}_i(X_1,\dots,X_n), \quad \forall i\in [n], \pi \in S_n. 
\end{align*}
In fact, any estimator $\widehat{\theta}$ can be symmetrized into a PI estimator $\widehat{\theta}^{\rm PI}(X^n) = \frac{1}{n!}\sum_{\pi\in S_n} \pi^{-1}\circ \widehat{\theta}(\pi\circ X^n)$ without increasing the worst-case MSE over the orbit $\{\pi\circ \theta: \pi\in S_n\}$. Hannan and Robbins \cite{hannan1955asymptotic} set out to evaluate EB estimators against the \emph{permutation-invariant oracle} $\thetahatPI$:
\begin{align}\label{eq:PI_oracle}
\MSE(\theta, \thetahatPI) \le \MSE(\theta, \widehat{\theta}), \quad \forall \widehat{\theta}\in \calD^{\rm PI}.
\end{align}
The PI oracle turns out to be the Bayes estimator in a postulated Bayes model; we defer discussions and equivalent definitions of the PI oracle to \Cref{subsec:EB_preliminary}. Since $\calD^{\rm S}\subseteq \calD^{\rm PI}$, a stronger asymptotic guarantee for the EB estimator would be
\begin{align}\label{eq:regret_against_PI}
\regPI(\widehat{\theta}^{\rm EB}) &:= \sup_{\theta\in \Theta} \pth{ \MSE(\theta,\widehat{\theta}^{\rm EB}) - \MSE(\theta,\thetahatPI) }= o(n).
\end{align}
This regret against the PI oracle can be viewed as more fundamental in compound settings, since its definition does not rely on EB approximations. Indeed, if \eqref{eq:regret_against_PI} holds, this shows that the bound \eqref{eq:PI_oracle} is asymptotically attainable uniformly in $\theta$ by a \emph{proper} and \emph{legal} rule, in the sense that the EB estimator $\widehat{\theta}^{\rm EB}$ belongs to $\calD^{\rm PI}$ and does not depend on the ground truth $\theta$. Both quantities $\regS$ and $\regPI$ are called \emph{compound regrets} in the literature, as opposed to the EB regret which assumes i.i.d. $\theta_1,\dots,\theta_n\sim G$. However, compared with the more widely used version in \eqref{eq:regret_against_simple}, far less is known about the larger quantity \eqref{eq:regret_against_PI}. Existing results \cite{hannan1955asymptotic,greenshtein2009asymptotic,han2024approximate} primarily establish the regret equivalence $\regS\asymp \regPI$ with a fixed support of $G_n$, and become ineffective in many interesting settings where the support of $G_n$ grows with $n$, even moderately.

At a high level, the target regret equivalence $\regS\asymp \regPI$ asks how well the simpler separable rule approximates the optimal PI rule. This question naturally reduces to a \emph{mean-field} approximation, as witnessed by the exact expressions of both oracles in a \emph{postulated} Bayes model, in which the parameter vector $\theta=(\theta_1,\dots,\theta_n)$ is randomly permuted before generating observations. Specifically, let $\pi\sim\Unif(S_n)$ be a random permutation, and define $\widetilde{\theta}_i = \theta_{\pi(i)}$, with independent observations $X_i\sim P_{\widetilde{\theta}_i}$. Under this model, \cite{greenshtein2009asymptotic} showed that
\begin{align}\label{eq:postulated_bayes_model}
    \thetahatS_i = \bE[\widetilde{\theta}_i|X_i], \qquad \thetahatPI_i = \bE[\widetilde{\theta}_i|X^n], \qquad \forall i\in [n].
\end{align}
Note that if $(\widetilde{\theta}_i, X_i)$ \emph{were} truly independent across $i\in [n]$, the two estimators would coincide. Their difference arises solely due to dependence induced by the permutation, and this difference is exactly the type of approximation error captured by permutation mixtures. Hence, understanding when $\thetahatPI \approx \thetahatS$, i.e. when mean-field approximations are tight in this postulated Bayes model, is essential for tightening the theoretical guarantees of EB estimators. This connection is the central motivation behind our refined analysis of permutation mixtures in the remainder of this paper.

\subsection{Notation}
Throughout the paper all logarithms are in base $e$. For a positive integer $n$, let $[n]:=\{1,\ldots,n\}$, and $S_n$ be the symmetric group over $[n]$. For a vector $x$, let $\|x\|$ be its $\ell_2$ norm. For a square matrix $A = (a_{ij})_{i,j\in [n]}\in \mathbb{R}^{n\times n}$, let $\trace(A) = \sum_{i=1}^n a_{ii}$ be its trace, and $\Perm(A) = \sum_{\pi \in S_n} \prod_{i=1}^n a_{i\pi(i)}$ be its permanent. For probability measures $P$ and $Q$ on the same probability space, let 
\begin{align*}
	&\TV(P,Q) = \frac{1}{2}\int |\rmd P - \rmd Q|, \qquad H^2(P,Q) = \int \pth{\sqrt{\rmd P} - \sqrt{\rmd Q}}^2\end{align*}
be the total variation (TV) and squared Hellinger distances, respectively, and 
\begin{align*}
	&\KL(P\|Q) = \int \rmd P\log \frac{\rmd P}{\rmd Q}, \qquad \chi^2(P\|Q) = \int \frac{(\rmd P - \rmd Q)^2}{\rmd Q}
\end{align*}
be the Kullback--Leibler (KL) and $\chi^2$ divergences, respectively. A collection of inequalities between the above distances/divergences can be found in \cite[Chapter 7.6]{polyanskiy2024information}. 

We shall use the following standard asymptotic notations. For non-negative sequences $\{a_n\}$ and $\{b_n\}$, let $a_n = O(b_n)$ denote $\limsup_{n\to\infty} a_n/b_n < \infty$, and $a_n = o(b_n)$ denote $\limsup_{n\to\infty} a_n/b_n = 0$. In addition, we write $a_n = \Omega(b_n)$ for $b_n = O(a_n)$, $a_n = \omega(b_n)$ for $b_n = o(a_n)$, and $a_n = \Theta(b_n)$ for both $a_n = O(b_n)$ and $b_n = O(a_n)$. We will also use the notations $O_\theta, o_\theta$, etc. to denote that the hidden factor depends on some external parameter $\theta$.

\subsection{Main results}
\paragraph{Dimension-independent bounds} Given a class $\calP$ of distributions, we define the best dimension-independent $\chi^2$ bound as
\begin{align}\label{eq:dimension_indep_bound}
\chi^2(\calP) = \sup_{n}\sup_{P_1,\dots,P_n\in \calP} \chi^2(\bP_n \| \bQ_n), 
\end{align}
where $\bP_n$ and $\bQ_n$ are defined in \cref{eq:bP_n} and \cref{eq:bQ_n} in \Cref{def:perm_mixture}, respectively, determined by $P_1,\dots,P_n$. To present our upper bound on $\chi^2(\calP)$, we introduce the following quantities. 

\begin{definition}[Quantities of $\calP$]\label{defn:quantity}
For a given family $\calP$ of probability distributions over the same space, define: 
\begin{enumerate}
    \item The \emph{$\chi^2$ channel capacity}, denoted by $\Capa$: 
    \begin{align}\label{eq:capacity}
    \Capa &:= \sup_{\rho} I_{\chi^2}(\theta; X) = \sup_{\rho} \bE_{\theta\sim \rho}\qth{\chi^2(P_{\theta} \| \bE_{\theta'\sim \rho}[P_{\theta'}]) }, 
    \end{align}
    where $(P_{\theta})_{\theta} = \calP$ is a parametrization of $\calP$, and $\rho$ is a prior distribution of $\theta$, with $X|\theta\sim P_\theta$; 
    \item The \emph{R\'{e}nyi partition diameter} of order $k\in \naturals$, denoted by $\Diam$: 
    \begin{align}\label{eq:partition-diameter}
    \Diam := \inf_{\text{\rm partitions } \calP = \bigsqcup_{i=1}^k \calP_i} \max_{i\in [k]}\sup_{P, Q\in \calP_i} D_{1/2}(P, Q), 
    \end{align}
    where $D_{1/2}(P,Q) = -2\log(1-\frac{H^2(P,Q)}{2})$ is the R\'{e}nyi divergence of order $\frac{1}{2}$.\footnote{Due to symmetry, we write $D_{1/2}(P,Q)$ instead of the classic notation $D_{1/2}(P\|Q)$. However, the readers should still be warned that $D_{1/2}(P,Q)$ is \emph{not} a metric, for the triangle inequality does not hold.} 
\end{enumerate}
\end{definition}

The notion of $\chi^2$ channel capacity was introduced in \cite{han2024approximate}, which replaces the KL divergence by the $\chi^2$ divergence in the usual definition of channel capacity for the channel class $\calP$ from $\theta$ to $X$. The notion of R\'{e}nyi partition diameter in \eqref{eq:partition-diameter} is new, and describes how the distribution class $\calP$ can be partitioned into $k$ subclasses, each of a small diameter under the R\'{e}nyi divergence $D_{1/2}$. Note that when $k=1$, 
\begin{align*}
e^{\sfD_1(\calP)} = \pth{1-\frac{\max_{P,Q\in \calP} H^2(P,Q)}{2}}^{-2}
\end{align*}
coincides with the definition of \emph{maximum $H^2$ singularity} in \cite{han2024approximate}; moreover, $\Diam \le \sfD_1(\calP)<\infty$ if no pairs of distributions in the compactification of $\calP$ under the Hellinger distance are mutually singular. Finally, we note that both quantities $\Capa$ and $\Diam$ depend only on $\calP$ but not on the dimension $n$ in the definitions of $\bP_n$ and $\bQ_n$. We present our first theorem below, which is our main upper bound on the dimension-independent quantity $\chi^2(\calP)$. 

\begin{thm}[Dimension-independent upper bound]\label{thm:main_upper}
The following upper bound holds:
\begin{align*}
&\log(1+\chi^2(\calP))  \le C \Big(\sum_{k=1}^{\lfloor \Capa \rfloor+1} \Diam + (\Capa+1) \log_+\log \Capa \Big),
\end{align*}
where $C>0$ is a universal constant, and $\log_+(x) := \log(x\vee e)$. 
\end{thm}


As a comparison, the best upper bound in \cite{han2024approximate} for the case $\Capa \ge 1$ is $\log(1 + \chi^2(\mathcal{P})) \le \Capa (\sfD_1(\mathcal{P}) + 1)$. Since the function $k \mapsto \sfD_k(\mathcal{P})$ is non-increasing, the first term in \Cref{thm:main_upper} is always at least as good as this existing bound, and offers improvements when $\sfD_k(\mathcal{P})$ decays rapidly with $k$. In particular, if $\sfD_k(\mathcal{P}) \asymp \sfD_1(\mathcal{P}) / k^{\alpha}$ for some $\alpha > 1$, then this sum remains $O(\sfD_1(\mathcal{P}))$ regardless of how large $\Capa$ is. The second term in \Cref{thm:main_upper} scales linearly with $\Capa$ (up to an additional $\log\log$ factor), but is typically dominated by the first term; see the examples in \Cref{cor:dim-indep-bound}.

{\blue Intuitively, \cref{thm:main_upper} can be interpreted in terms of an integral of a covering number. For $r\ge 0$, let $N(r):= \min\{k\geq 1: \Diam \le r\}$ denote the smallest integer $k$ such that $\calP$ can be partitioned into $k$ subclasses, each having diameter at most $r$ under the R\'enyi divergence $D_{1/2}$. Thus, $N(r)$ may be viewed as a covering number of $\calP$ at radius $r$ under $D_{1/2}$. For many classes $\calP$ with $\Capa\ge 1$, including all the examples in \Cref{cor:dim-indep-bound} below, the $\chi^2$ channel capacity $\Capa$ is of the same order as the covering number $N(1)$ at radius $1$. Since
\begin{align*}
 \sum_{k=1}^{N(1)}\Diam
&= N(1) + \sum_{k=1}^{N(1)} \int_{\sfD_{N(1)}(\calP)}^{\sfD_1(\calP)} |\sth{k\ge 1: \sfD_k(\calP) > r}| \rmd r \\&\asymp N(1)+ \int_1^{\sfD_1(\calP)} N(r)\rmd r, 
\end{align*}
the upper bound in \Cref{thm:main_upper} is roughly
\begin{align}\label{eq:entropy-integral}
\log(1+\chi^2(\calP)) \lesssim \int_1^{\sfD_1(\calP)} N(r)\rmd r + N(1)\log_+\log N(1),
\end{align}
which involves integrating the covering number $N(r)$ from $r=1$ to $r=\sfD_1(\calP)$. By contrast, the upper bound in \cite{han2024approximate} can be interpreted as $\int_1^{\sfD_1(\calP)} N(1)\rmd r + N(1)$. Thus, the improvement of \Cref{thm:main_upper} over \cite{han2024approximate} may be viewed intuitively as the improvement achieved by a chaining argument over a one-stage covering argument. As a concrete example, consider the normal mean model $\calP = \{ \calN(\theta,1): |\theta|\le \mu \}$ with $\mu>1$. In this case, $\Capa\asymp \mu $ and $\Diam\asymp \frac{\mu^2}{k^2}$. Therefore, $N(r)\asymp \frac{\mu}{\sqrt{r}}$, and indeed $N(1)\asymp \Capa$. Moreover, $ \int_1^{\mu^2} N(r)\rmd r\asymp \mu^2$, whereas $\int_1^{\mu^2} N(1)\rmd r\asymp \mu^3$.
}

To show the tightness of \Cref{thm:main_upper} in many examples, we also derive a general lower bound of the quantity $\chi^2(\calP)$. We recall the definition of the channel overlap matrix in \cite{han2024approximate, kunisky2024low}. 
\begin{definition}[Channel overlap matrix]\label{defn:channel-overlap-matrix}
Given $n$ probability distributions $P_1,\dots,P_n$ on the same space, the channel overlap matrix $A\in \bR^{n\times n}$ is defined by
\begin{align*}
A_{ij} = \frac{1}{n}\int \frac{\rmd P_i \rmd P_j}{\rmd \overline{P}}, \quad i,j\in [n], 
\end{align*}
where $\overline{P} := \frac{1}{n}\sum_{i=1}^n P_i$ is the average distribution. 
\end{definition}

In \cite[Lemma 5.2]{han2024approximate}, it is shown that $A$ is doubly stochastic and PSD, with an eigenstructure $1=\lambda_1(A)\ge \lambda_2(A)\ge \cdots \ge \lambda_n(A)\ge 0$. In addition, the $\chi^2$ divergence is related to the permanent of $A$ via the identity \cite[Lemma 5.1]{han2024approximate}
\begin{align}\label{eq:permanent}
    1 + \chi^2(\bP_n \| \bQ_n) = \frac{n^n}{n!}\Perm(A). 
\end{align}
The next theorem presents a general lower bound of $\chi^2(\calP)$. 

\begin{thm}[Dimension-independent lower bound]\label{thm:main_lower}
The following lower bound holds:
\begin{align*}
1 + \chi^2(\calP) \ge \sup_n \sup_{P_1,\dots,P_n\in \calP} \prod_{i=2}^n \frac{1}{\sqrt{1-\lambda_i(A)^2}}, 
\end{align*}
where $A=A(P_1,\dots,P_n)$ is the channel overlap matrix in \Cref{defn:channel-overlap-matrix}. In particular, 
\begin{align*}
1 + \chi^2(\calP) \ge \sup_n \sup_{P_1,\dots,P_n\in \calP}  \frac{1}{\sqrt{n}}\prod_{i=2}^n \frac{1}{\sqrt{1-A_{ii}^2}}. 
\end{align*}
\end{thm}

We compare \Cref{thm:main_lower} with two results in \cite{han2024approximate} that involve the channel overlap matrix. First, a similar lower bound in \cite[Lemma 6.1]{han2024approximate} reads that
\begin{align}\label{eq:A-lower-old}
1+\chi^2(\calP) \ge \sup_n \sup_{P_1,\dots,P_n\in \calP} \frac{1}{\sqrt{1-\lambda_2(A)^2}}, 
\end{align}
which only involves the second largest eigenvalue of $A$. By contrast, the lower bound in \Cref{thm:main_lower} involves \emph{all} non-leading eigenvalues. Technically, \eqref{eq:A-lower-old} is proven by checking the variance of a linear test function, while our proof for \Cref{thm:main_lower} applies an involved saddle-point analysis to a limiting matrix permanent. This finding also parallels with an intermediate upper bound in \cite[Sec. 5.2]{han2024approximate}: using auxiliary complex normal random variables, it was established that
\begin{align}\label{eq:A-upper}
1 + \chi^2(\calP) \le \sup_n \sup_{P_1,\dots,P_n\in \calP} \prod_{i=2}^n \frac{1}{1-\lambda_i(A)}. 
\end{align}
Combining both upper and lower bounds, we see that controlling the gaps between all non-leading eigenvalues of $A$ and the leading eigenvalue $\lambda_1(A) = 1$ plays a significant role in bounding $\chi^2(\calP)$. In fact, the main technical contribution of \Cref{thm:main_upper} is to connect these gaps to the R\'{e}nyi partition diameters of $\calP$. 

Finally, by combining the upper and lower bounds, we obtain (near-)tight characterizations of $\chi^2(\mathcal{P})$ for a variety of distribution classes $\mathcal{P}$. Since \cite{han2024approximate} has already established the tightness of the relation $\chi^2(\mathcal{P}) = \Theta(\Capa^2)$ in the regime $\Capa \le 1$, the following examples focus exclusively on parameter regimes where $\Capa > 1$.

\begin{corollary}\label{cor:dim-indep-bound}
The following results hold for specific families:
\begin{enumerate}
    \item Normal mean family $\calP=\{\calN(\theta,1): |\theta|\le \mu\}$ with $\mu\ge 1$: 
    $\chi^2(\calP) = \exp(\Theta(\mu^2)).$ 
    \item Poisson family $\calP = \{\Poi(\lambda): m\le \lambda\le M\}$ with $\sqrt{M}-\sqrt{m}\ge 1$: $$
    \chi^2(\calP) = \exp(\Theta((\sqrt{M}-\sqrt{m})^2)). $$
    \item Normal variance family $\calP=\{\calN(0,\sigma^2): |\log \sigma|\le M \}$ with $M\ge 1$: 
    \begin{align*}
    \chi^2(\calP) \in \exp(\Omega(M) \cap O(M\log M)). 
    \end{align*}
    \item Normal mean family $\calP=\{\calN(\theta,I_d): \theta\in \bR^d, \|\theta\|\le \mu\}$ in $d\ge 2$ dimensions, with $\mu\ge 1$:
    \begin{align*}
    \chi^2(\calP) \in \begin{cases}
        \exp(\Omega(\mu^2)\cap O(\mu^2\log \mu)) &\text{if } d = 2, \\
        \exp(\Omega_d(\mu^d)\cap O_d(\mu^d\log\log \mu)) & \text{if } d \ge 3. 
    \end{cases}
    \end{align*}
    \item Discrete distribution family $\calP = \{p\in \Delta_m: p_1\ge \varepsilon\}$ with $m\ge 2$ and overlap $\varepsilon\in (0, \frac{1}{2})$: 
    \begin{align*}
    \chi^2(\calP) = \exp\pth{ \Theta\pth{m\log\frac{1}{\varepsilon}} }. 
    \end{align*}
\end{enumerate}
\end{corollary}

We note that, for all the above examples (except for Example 5 with a special geometry), the existing upper bounds in \cite{han2024approximate} fail to capture the correct polynomial scaling in $\mu$ or $M$ in the exponent. These examples also illustrate how different structures of information geometry affect the behavior of $\chi^2(\mathcal{P})$. For instance, in the normal mean family in $d\ge 1$ dimensions, one has $\Capa=O(\mu^d)$ and $\sfD_k(\mathcal{P}) \asymp \mu^2 / k^{2/d}$, so that the sum in \Cref{thm:main_upper} is $\widetilde{O}(\mu^{2 \vee d})$, revealing a phase transition between the regimes $d=1,2$ and $d \ge 3$. Moreover, applying Theorem \ref{thm:main_lower} yields a matching exponent, confirming the tightness of this transition. As another example, in Example 5, the R\'enyi partition diameter scales as $\sfD_k(\mathcal{P}) \asymp \log(1/\varepsilon)$ for all $1 \le k \le m-1$, i.e. it does not decay with $k$. This explains why the bound in \cite{han2024approximate}, which relies solely on $\sfD_1(\mathcal{P})$, remains tight in that case. In summary, these examples suggest that the behavior of $\chi^2(\mathcal{P})$ crucially depends on how the R\'enyi partition diameter decays with $k$. 

\paragraph{Refinements in finite dimensions} Similar to \eqref{eq:dimension_indep_bound}, we can define its dimension-dependent counterpart. For fixed $n\ge 2$, let
\begin{align}\label{eq:dimension_dep_bound}
\chi^2(\calP, n) = \max_{m\le n}\sup_{P_1,\dots,P_m\in \calP} \chi^2(\bP_m \| \bQ_m). 
\end{align}
As discussed in the introduction, the main advantage of studying the dimension-dependent quantity \eqref{eq:dimension_dep_bound} is to achieve a tighter dependence on $\mathcal{P}$ at the cost of a slow growth in $n$. The following theorem demonstrates that this is possible.

{\blue
\begin{thm}[Dimension-dependent upper bound]\label{thm:main_upper_dimdep}
The following upper bound holds: 
\begin{align*}
\log(1+\chi^2(\calP,n)) & \le C \Big(\sum_{k=1}^{\lfloor \Capa \rfloor+1} \min\{\Diam, \log n\} + (\Capa+1) \log_+\log \Capa \Big),
\end{align*}
where $C>0$ is a universal constant, and $\log_+(x) := \log(x\vee e)$. 
\end{thm}

Clearly, \Cref{thm:main_upper_dimdep} is an improvement of \Cref{thm:main_upper} by truncating the R\'enyi partition diameter $\Diam$ at a dimension-dependent threshold $\log n$. Using the integral view of \eqref{eq:entropy-integral}, \Cref{thm:main_upper_dimdep} roughly shows that
\begin{align}\label{eq:entropy-integral-finite-n}
\log(1+\chi^2(\calP,n)) \lesssim \int_1^{\min\{\sfD_1(\calP),\log n\}} N(r)\rmd r + N(1)\log_+\log N(1),
\end{align}
where compared with \eqref{eq:entropy-integral}, the upper limit of the integral improves from $r=\sfD_1(\calP)$ to $r=\min\{\sfD_1(\calP), \log n\}$. In the normal mean example with $N(r)\asymp \frac{\mu}{\sqrt{r}}$ and $\sfD_1(\calP) \asymp \mu^2$, the new integral in \eqref{eq:entropy-integral-finite-n} gives $O(\min\{\mu^2,\mu\sqrt{\log n}\})$, which exhibits a linear dependence on $\mu$ for small $n$. Although we do not have a general lower bound similar to \Cref{thm:main_lower} due to difficulty of handling a finite $n$, for both normal and Poisson models, \Cref{thm:main_upper_dimdep} turns out to yield tight bounds for all regimes of $n$. 
}

\begin{corollary}\label{cor:dim-dep-bound}
Let $n\ge 2$. The following results hold: 
\begin{enumerate}
    \item Normal mean family $\calP = \{\calN(\theta,1): |\theta|\le \mu\}$: 
    \begin{align*}
    \chi^2(\calP,n) = \begin{cases}
        \Theta(\mu^4) &\!\!\!\!\! \text{if } \mu \le 1, \\
        \exp(\Theta(\mu^2)) &\!\!\!\!\! \text{if } 1<\mu\le \sqrt{\log n}, \\
        \exp(\Theta(\mu\sqrt{\log n})) &\!\!\!\!\! \text{if } \sqrt{\log n} < \mu \le n^{0.99}. 
    \end{cases}
    \end{align*}
    \item Poisson family $\calP = \{\Poi(\lambda): 0\le \lambda\le M\}$: 
    \begin{align*}
    \chi^2(\calP,n) = \begin{cases}
        \Theta(M^2) &\!\!\!\!\! \text{if } M \le 1, \\
        \exp(\Theta(M)) &\!\!\!\!\! \text{if } 1<M\le \log n, \\
        \exp(\Theta(\sqrt{M\log n})) &\!\!\!\!\! \text{if } \log n < M \le n^{1.99}. 
    \end{cases}
    \end{align*}
\end{enumerate}
\end{corollary}

Corollary~\ref{cor:dim-dep-bound} reveals an additional phase transition in the behavior of the dimension-dependent quantity $\chi^2(\mathcal{P}, n)$. For the normal mean family, two distinct elbows occur at $\mu \asymp 1$ and $\mu \asymp \sqrt{\log n}$, with the exponent shifting from $\mu^2$ to $\mu$ at the second elbow. A similar transition is observed for the Poisson family. This highlights the importance of considering the dimension-dependent quantity $\chi^2(\mathcal{P}, n)$ in applications where the dimension $n$ is moderate: the regime identified in Corollary~\ref{cor:dim-dep-bound} may indeed be the most relevant. For instance, this new regime will play a key role in our mean-field analysis for compound decision problems. 

The upper bounds in Corollary~\ref{cor:dim-dep-bound} follow from \cite{han2024approximate} when $\mu\le 1$ or $M\le 1$, and Theorem~\ref{thm:main_upper_dimdep} otherwise. Lower bounds are obtained via carefully constructed test functions and are deferred to the appendix. Finally, the condition $\mu \le n^{0.99}$ (or $M \le n^{1.99}$) is not superfluous: even for extremely large $\mu$ (or $M$), a deterministic upper bound
\begin{align*}
\chi^2(\calP,n) \le \frac{n^n}{n!} - 1 = e^{O(n)}
\end{align*}
always holds, so the exponent cannot exceed $O(n)$. {\blue The choices of the constants $0.99$ and $1.99$ are not important, and can be replaced by any fixed exponents $a<1$ and $b<2$, respectively.}

\paragraph{Applications to compound decision problems} The improvement in the dimension-independent bound of \Cref{thm:main_upper}, along with the identification of a new dimension-dependent regime in \Cref{cor:dim-dep-bound}, turns out to be critical for establishing tight regret guarantees in compound decision problems. Recall the general setup: let ${P_\theta : \theta \in \Theta \subseteq \mathbb{R}}$ be a parametric family of distributions. The learner observes $X_i \sim P_{\theta_i}$ independently for $i \in [n]$ and aims to estimate the unknown parameters $\theta_1, \dots, \theta_n$. For an estimator $\widehat{\theta}$, we define its regret relative to the best separable and permutation-invariant (PI) decision rules as
\begin{align}\label{eq:regret_S}
\regS(\widehat{\theta}) &:= \sup_{\theta\in\Theta} \pth{ \MSE(\theta,\widehat{\theta}) - \MSE(\theta,\thetahatS) } \\ \label{eq:regret_PI}
\regPI(\widehat{\theta}) &:= \sup_{\theta\in\Theta} \pth{ \MSE(\theta,\widehat{\theta}) - \MSE(\theta,\thetahatPI) }
\end{align}
respectively, where $\MSE$ stands for the mean squared error in \eqref{eq:MSE}, and $\thetahatS, \thetahatPI$ are the separable and PI oracles defined in \eqref{eq:separable_oracle} and \eqref{eq:PI_oracle}, respectively. While upper bounds on $\regS$ have been the focus of substantial work on empirical Bayes methods \cite{Jiang_2009,brown2009nonparametric,polyanskiy2021sharp,shen2022empirical,jana2022optimal,jana2023empirical,ghosh2025stein}, much less is known about the other compound regret $\regPI$. 

In this paper, we bridge this gap by directly bounding the difference $\regPI(\widehat{\theta}) - \regS(\widehat{\theta})$, thereby showing that these two notions of compound regret are almost equivalent. We focus on Gaussian and Poisson models, as they serve as prototypical examples in the literature.

In the Gaussian location model, we have independent $X_i\sim \calN(\theta_i,1)$ for $i\in [n]$. The following theorem summarizes the upper bound on the regret difference with various tail conditions on $\theta$.
\begin{thm}[Gaussian regret transfer]\label{thm:EB-gaussian}
In the Gaussian location model, for every estimator $\widehat{\theta}$, it holds that
\begin{align*}
0 \le \regPI(\widehat{\theta}) - \regS(\widehat{\theta}) \le r_n, 
\end{align*} 
where: 
\begin{enumerate}
    \item $r_n = O(h^4  \wedge h\log^{3/2} n)$ if $\theta_i \in [-h,h]$ for all $i\in [n]$;
    \item $r_n = O_s(\log^2 n)$ if the empirical distribution $G_n=\frac{1}{n}\sum_{i=1}^n \delta_{\theta_i}$ is $s$-subGaussian; 
    \item $r_n = O_p(n^{\frac{1}{1+p}}\log^{\frac{4+3p}{2(1+p)}} n)$ if $G_n$ lies in the weak $\ell_p$ ball with some $p>0$, i.e. 
    \begin{align}\label{eq:weak-ell-p-ball}
    \mu_p^w(G_n) := \bpth{\sup_{x>0} x^p \int_{|u|>x} G_n(\rmd u) }^{1/p} = O(1). 
    \end{align}
\end{enumerate}
\end{thm}

Under the same tail conditions, the following bounds are known for {\blue the minimax regret $\regS_\star:=\inf_{\widehat{\theta}}\regS(\widehat{\theta})$}:
\begin{enumerate}
    \item If $\theta_i\in [-h,h]$ for all $i\in [n]$, the minimax lower bound is $\regS_\star=\Omega_h((\frac{\log n}{\log\log n})^2)$ \cite{polyanskiy2021sharp}; a recent estimator attains this upper bound under a different EB setting \cite{chen2026sharp}. 
    \item If $G_n$ is $s$-subGaussian, the best known upper bound is $\regS_\star = O_s(\log^5 n)$ \cite{Jiang_2009} (or an upper bound of $O_s((\log n)^3(\log\log n)^2)$ for $n$ times the EB regret \cite{chen2026sharp}), and the best known minimax lower bound is $\regS_\star = \Omega_s(\log^2 n)$ \cite{polyanskiy2021sharp}.
    \item If $G_n$ belongs to the weak $\ell_p$ ball, the best known upper bound is $\regS_\star=O_p(n^{\frac{1}{1+p}}\log^{\frac{8+9p}{2+2p}} n)$ \cite{Jiang_2009}, and we show a minimax lower bound $\regS_\star = \Omega_p(n^{\frac{1}{1+p}}\log^{-\frac{p}{2(1+p)}} n)$ in Appendix \ref{append:regret_LB}. 
\end{enumerate}
Combining these known results and \Cref{thm:EB-gaussian}, we have the following corollary. 
\begin{corollary}[Gaussian regret equivalence]\label{cor:EB-gaussian}
Under the setting of \Cref{thm:EB-gaussian}, for every estimator $\widehat{\theta}$ we have:
\begin{itemize}
    \item $\regPI(\widehat{\theta})\asymp \regS(\widehat{\theta})$ in Cases 1 and 2; 
    \item $\regPI(\widehat{\theta})=\widetilde{\Theta}( \regS(\widehat{\theta}))$ in Case 3. 
\end{itemize}
\end{corollary}

Similarly, we also have the following guarantees for the regret difference in the Poisson model, where $X_i\sim \Poi(\theta_i)$ are independent. 

\begin{thm}[Poisson regret transfer]\label{thm:EB-Poisson}
In the Poisson model, {\blue for every estimator $\widehat{\theta}$,} it holds that
\begin{align*}
0 \le \regPI(\widehat{\theta}) - \regS(\widehat{\theta}) \le r_n, 
\end{align*} 
where: 
\begin{enumerate}
    \item $r_n = O(h^3  \wedge (h\log n)^{3/2})$ if $\theta_i \in [0,h]$ for all $i\in [n]$;
    \item $r_n = O_s(\log^3 n)$ if the empirical distribution $G_n=\frac{1}{n}\sum_{i=1}^n \delta_{\theta_i}$ is $s$-subexponential; 
    \item $r_n = O_p(n^{\frac{3}{2p+1}}\log^{\frac{3(p+1)}{2p+1}} n)$ if $G_n$ has bounded $p$-th moment for some $p\ge 1$, i.e. $\int u^pG_n(\rmd u) = O(1)$.
\end{enumerate}
\end{thm}
The following bounds are known for the {\blue minimax regret $\regS_\star:=\inf_{\widehat{\theta}}\regS(\widehat{\theta})$} in the Poisson model: 
\begin{enumerate}
    \item If $\theta_i \in [0,h]$ for all $i\in [n]$, the minimax lower bound is $\regS_\star = \Omega_h((\frac{\log n}{\log\log n})^2)$ \cite{polyanskiy2021sharp}; Robbins's estimator attains this bound under a different EB setting. 
    \item If $G_n$ is $s$-subexponential, the minimax lower bound is $\regS_\star = \Omega_s(\log^3 n)$ \cite{polyanskiy2021sharp}; Robbins's estimator attains this bound under a different EB setting. 
    \item If $G_n$ has bounded $p$-th moment, the best known upper and lower bounds are $O_p(n^{\frac{3}{2p+1}}\log^{13}n)$ and $\Omega_p(n^{\frac{3}{2p+1}}\log^{-11}n)$, respectively\footnote{Technically, the lower bound in \cite{shen2022empirical} was established for the EB regret. However, it is straightforward to adapt their proof to the compound regret.} \cite{shen2022empirical}.
\end{enumerate}
Combining these known results and \Cref{thm:EB-Poisson}, we have the following corollary. 
\begin{corollary}[Poisson regret equivalence]\label{cor:EB-Poisson}
Under the setting of \Cref{thm:EB-Poisson}, {\blue for every estimator $\widehat{\theta}$,} we have:
\begin{itemize}
    \item $\regPI(\widehat{\theta})\asymp \regS(\widehat{\theta})$ in Cases 1 and 2; 
    \item $\regPI(\widehat{\theta})=\widetilde{\Theta}( \regS(\widehat{\theta}))$ in Case 3. 
\end{itemize}
\end{corollary}

In summary, \Cref{cor:EB-gaussian} and \Cref{cor:EB-Poisson} establish strong equivalences between two notions of compound regret in both the Gaussian and Poisson models, thereby validating the use of EB-style approximations for the PI oracle in compound decision problems. 

\subsection{Related work}\label{subsec:related-work}
Approximating exchangeable distributions by (mixtures of) product measures is a classical heuristic dating back to de Finetti \cite{de1929funzione} and Hewitt--Savage \cite{hewitt1955symmetric}. When $P_1, \dots, P_n$ are point masses, comparing $\bP_n$ and $\bQ_n$ reduces to comparing sampling without replacement and sampling with replacement, a well-studied topic in works by Stam \cite{stam1978distance} and Diaconis--Freedman \cite{diaconis1980finite,diaconis1987dozen}, with more recent developments in \cite{gavalakis2021information,johnson2024relative}. In this case, contiguity results hold only for $k$-dimensional marginals of $\bP_n$ and $\bQ_n$, with $k = O(\sqrt{n})$ or $k \le (1 - \varepsilon)n$ if the support is bounded.

A surprising recent development was the discovery that, in certain “noisy” settings, the full joint distribution satisfies the contiguity relation $\{\bP_n\} \triangleleft \{\bQ_n\}$. This was first shown in \cite{Din22} (Ph.D. thesis) for a balanced two-point Gaussian/Bernoulli scenario, and independently in \cite{tang2023capacity} when $\{P_1, \dots, P_n\}$ consists of a bounded number of discrete distributions, motivated by converse bounds for the capacity of noisy permutation channels. The former used a convexity argument tailored to the toy case, while the latter employed method-of-types and anticoncentration for Poisson–binomial distributions, which unavoidably incurs a linear dependence on the cardinality of $\{P_1, \dots, P_n\}$. The recent work \cite{han2024approximate} generalized these results substantially, by establishing $\chi^2(\bP_n\|\bQ_n) = O(1)$ in much broader settings (including when $P_1,\dots,P_n$ are continuous and all distinct), via two functional approaches. Our work builds on their functional inequality \eqref{eq:A-upper}, but identifies the higher-order Cheeger inequality as an appropriate tool to obtain sharper dimension-independent bounds (\Cref{thm:main_upper} and \Cref{cor:dim-indep-bound}), and complements it with a lower bound of similar form (\Cref{thm:main_lower}) obtained through a saddle-point analysis of a limiting matrix permanent.

These refined $\chi^2$ bounds are not purely of information-theoretic interest: they play a critical role in compound decision problems, where tight divergence control for permutation mixtures enables sharper regret bounds. In particular, \cite{han2025best} was the first to refine $\KL(\bP_n \| \bQ_n)$ in moderate dimensions and apply it to $\regPI$. {\blue The technique used in \cite{han2025best} and in an earlier version of the present paper \cite{LH25} partitions $\calP$ into $k$ subclasses, applies a method-of-types calculation, and then optimizes over $k$. By contrast, the current \Cref{thm:main_upper_dimdep} provides a more direct analogue of the dimension-independent result in \Cref{thm:main_upper} and avoids this auxiliary optimization step.} However, due to the absence of tight dimension-independent bounds, their results do not capture the optimal dependence on the $\log n$ factor, nor do they identify the tightness of this approach or the sharp phase transitions in the Gaussian and Poisson families (\Cref{cor:dim-dep-bound}). As a result, their regret bounds had loose logarithmic factors, i.e. $O(\log^{2.5} n)$ in the Gaussian model with subGaussian tails and $O(\log^4 n)$ in the Poisson model with subexponential tails, which are insufficient to establish strong regret equivalence $\regPI \asymp \regS$ in \Cref{cor:EB-gaussian} and \Cref{cor:EB-Poisson}. In comparison, our mean-field analysis bounds the $\chi^2$ divergence (which upper bounds KL) without relying on the chain rule (\Cref{thm:main_upper_dimdep}), and produces tight bounds that are essential for achieving the sharp results in compound decision problems.


In statistics, empirical Bayes has long been the dominant framework for studying compound decision problems, dating back to \cite{robbins1951asymptotically,hannan1955asymptotic,robbins1956empirical}; we refer to the review articles \cite{casella1985introduction,zhang2003compound,efron2024empirical} for comprehensive overviews. A major theoretical breakthrough by \cite{Jiang_2009} provided near-optimal non-asymptotic bounds for $\regS$ in Gaussian models, inspiring a flurry of subsequent studies in Gaussian \cite{brown2009nonparametric,saha2020nonparametric,polyanskiy2021sharp,ghosh2025stein}, Poisson \cite{brown2013poisson,polyanskiy2021sharp,jana2022optimal,shen2022empirical,jana2023empirical}, and multiple testing contexts \cite{heller2021optimal,basu2021empirical,nie2023large,ignatiadis2025empirical}. While numerous EB estimators have been proposed, including the nonparametric MLE \cite{Jiang_2009}, Robbins's estimator \cite{robbins1956empirical}, Good--Turing estimator \cite{good1953population}, and ERM approaches \cite{jana2023empirical,ghosh2025stein}, all regret bounds in this literature target $\regS$ against the separable oracle. This leaves the more compelling benchmark $\regPI$ largely unexplored, mostly due to the technical complication of the PI oracle. 

Approximating the PI oracle by the separable oracle has a parallel line of study. The earliest asymptotic equivalence $\regPI - \regS = o(n)$ was shown in \cite{hannan1955asymptotic} via de Finetti arguments. The first non-asymptotic bound $\regPI - \regS = O(1)$ was obtained in \cite{greenshtein2009asymptotic} using a clever convexity trick (see a clearer exposition in \cite[Sec. 2.3]{han2024approximate}); however, this bound crucially requires that the problem parameters do \emph{not} grow with $n$. For example, in case 1 of \Cref{thm:EB-gaussian}, this $O(1)$ factor becomes $e^{O(h^2)}$, which is prohibitively large when $h \gg \sqrt{\log n}$. Using permutation mixtures, \cite{han2024approximate} obtained an upper bound with the same $e^{O(h^2)}$ scaling in the Gaussian model. Other related results include \cite{orlitsky2015competitive} and \cite{nie2023large}. The former circumvents the difficulty of analyzing the PI oracle by introducing a stronger ``natural oracle'', but this approach is limited to the Poisson model and does not yield tight regret rates. In a different context, \cite[Lemma~S9]{nie2023large} shows that $\regPI - \regS = \widetilde{O}(|\{\theta_1,\dots,\theta_n\}|)$ via method-of-types arguments, but this bound can be as large as $O(n)$ when $\theta_1,\dots,\theta_n$ are all distinct. 

A recent breakthrough by \cite{han2025best} used a novel combination of an interpolation step between each $\theta_i$ and $X_i$, the information-theoretic arguments of \cite{nie2023large}, and the permutation mixture framework of \cite{han2024approximate}. This led to a tight characterization of $\regPI$ in a Poisson model with KL loss, closing the gap left by \cite{orlitsky2015competitive}. While our analysis of $\regPI$ draws on several ideas from \cite{han2025best}, we leverage our sharper mean-field approximations for permutation mixtures to obtain tighter mean-field approximations for PI oracles. Crucially, these yield small regret gaps, i.e. $O(\log^2 n)$ in the Gaussian model under subGaussian tails and $O(\log^3 n)$ in the Poisson model under subexponential tails, that are necessary for establishing the strong regret equivalence $\regPI\asymp \regS$. Moreover, our argument unifies the bounds in \cite{greenshtein2009asymptotic} and \cite{han2025best}. In Case~1 of \Cref{thm:EB-gaussian}, we simultaneously recover the $O(h^4)$ bound for small $h = O(\sqrt{\log n})$ (where \cite{greenshtein2009asymptotic} proved it only for $h \le 1$) and improve the bound for large $h$ to $O(h \log^{3/2} n)$ (sharpening the logarithmic factor from \cite{han2025best}). Similar improvements hold in Poisson models: even for $h \le 1$, the approach of \cite{greenshtein2009asymptotic} cannot establish our $O(h^3)$ bound, due to an unbounded $\chi^2$ divergence between $\Poi(h)$ and $\Poi(0)$.

\subsection{Organization}
The remainder of this paper is organized as follows. In \Cref{sec:bounds}, we develop our improved bounds for permutation mixtures, by using the higher-order Cheeger inequality to prove the dimension-independent upper bound in \Cref{thm:main_upper}, a covering and method-of-types argument to establish the dimension-dependent upper bound in \Cref{thm:main_upper_dimdep}, and a saddle-point analysis to derive the lower bound in \Cref{thm:main_lower}. In \Cref{sec:EB}, we review the definitions of the separable and PI oracles, and show how the refined analysis of permutation mixtures leads to the regret guarantees in \Cref{thm:EB-gaussian} and \Cref{thm:EB-Poisson}. Open questions are discussed in \Cref{sec:discussions}, and the proofs of corollaries and auxiliary lemmas are deferred to the appendices.

%% file: Improved_Bounds_for_Permutation_Mixtures.tex
\section{Improved bounds for permutation mixtures}\label{sec:bounds}
In this section, we establish the divergence upper bounds in \Cref{thm:main_upper} and \Cref{thm:main_upper_dimdep}, and the lower bound in \Cref{thm:main_lower}. Both upper bounds depend on the information geometry of $\calP$, which is connected to the spectrum of the channel overlap matrix $A$ via a higher-order Cheeger inequality. By contrast, the lower bound is more analytical and obtained through a saddle-point analysis of a limiting matrix permanent. 

\subsection{Dimension-independent upper bounds via higher-order Cheeger inequality}
In this section we prove \Cref{thm:main_upper}, where a key technical result is the following lower bound on the eigen-gaps for all eigenvalues of the channel overlap matrix $A$. 
\begin{lemma}\label{lemma:spectral_gap}
  For $P_1,\dots,P_n\in \calP$, let A $\in \mathbb{R}^{n \times n}$ be the channel overlap matrix in \Cref{defn:channel-overlap-matrix}. Then for all $k\ge 1$, 
  \begin{align*}
    1 - \lambda_{10k}(A) \geq c\frac{e^{-2\Diam}} {\log (5k)}, 
  \end{align*}
  where $c>0$ is a universal constant. 
\end{lemma}

We show that \Cref{lemma:spectral_gap} implies \Cref{thm:main_upper}. First, by \eqref{eq:A-upper}, we have
\begin{align}\label{eq:objective}
\log(1+\chi^2(\calP)) \le \sup_{n}\sup_{P_1,\dots,P_n\in \calP} \sum_{i=2}^n \log \frac{1}{1-\lambda_i(A)}, 
\end{align}
where $1\ge \lambda_2(A)\ge \dots\ge \lambda_n(A)\ge 0$ are the eigenvalues of $A$ \cite[Lemma 5.2]{han2024approximate}. {\blue For $2\le i\le 9$, \cite[Lemma 5.2]{han2024approximate} shows that $1-\lambda_i(A)\ge e^{-\sfD_1(\calP)}$, so that $
\sum_{i=2}^9 \log \frac{1}{1-\lambda_i(A)} \le 8\sfD_1(\calP). $
For $i\ge 10$,} let $g_k := \max_{\ell\le k} ce^{-2\sfD_{\ell}(\calP)}/\log(5\ell)$, then \Cref{lemma:spectral_gap} tells that
\begin{align}\label{eq:constraint_1}
0\le \lambda_{10k}(A), \dots, \lambda_{10k+9}(A) \le 1-g_k, \quad \forall k\ge 1. 
\end{align}
By \cite[Lemma 5.2]{han2024approximate}, the trace of $A$ satisfies an additional constraint: 
\begin{align}\label{eq:constraint_2}
\sum_{i=10}^n \lambda_i(A) \le \trace(A) - 1 \le \Capa. 
\end{align}
Next we maximize the RHS of \eqref{eq:objective} over the variables $(x_{10},\ldots,x_n)=
(\lambda_{10}(A),\ldots,\lambda_n(A))
\in \mathbb{R}^{n-9}$, subject to the constraints \eqref{eq:constraint_1} and \eqref{eq:constraint_2}. Since $x\mapsto \log\frac{1}{1-x}$ is increasing and convex on $[0,1)$, the objective value increases under both operations $(x_i, x_j) \to (x_i+\varepsilon, x_j-\varepsilon)$ for $x_i\ge x_j$ and $\varepsilon > 0$, and $x_i\to x_i+\varepsilon$. Therefore, no such operations can be performed for the maximizer $(x_{10}^\star,\dots,x_n^\star)$, which exists by compactness of $[0,1-g_1]^{n-9}$. As $g_k\le g_{k+1}$, this maximizer must take the form
\begin{align*}
x_{10k}^\star, \dots, x_{10k+9}^\star \begin{cases}
   = 1-g_k & \text{for } k < k_0 \\
   \in [0,1-g_k] & \text{for } k = k_0 \\
   = 0 & \text{for } k > k_0
\end{cases}, \qquad \text{with }\sum_{i=10}^n x_i^\star = \Capa. 
\end{align*}
Since $g_k\le \frac{1}{2}$ (by choosing a small constant $c$ in \Cref{lemma:spectral_gap} if necessary), the cutoff $k_0$ satisfies
\begin{align*}
\Capa \ge \sum_{k=1}^{k_0-1} 10(1-g_k) \ge k_0-1 \Longrightarrow k_0 \le 1 + \Capa. 
\end{align*}
Therefore, the maximum value of the objective in \eqref{eq:objective} is at most
\begin{align*}
\sum_{i=10}^n \log \frac{1}{1-x_i^\star} \le \sum_{k=1}^{k_0} 10\log \frac{1}{g_k}  &\le 10\sum_{k=1}^{k_0} \pth{ 2\Diam + \log\frac{\log(5k)}{c}} \\
&= O\bpth{ \sum_{k=1}^{\lfloor \Capa\rfloor + 1} \Diam + (\Capa+1)\log_+\log \Capa }, 
\end{align*}
which is the claimed result of \Cref{thm:main_upper}. 

In the remainder of this section we prove \Cref{lemma:spectral_gap}, which lower bounds an algebraic quantity $1-\lambda_{2k}(A)$ by a function of the geometric quantity $\Diam$. To this end, we view the non-negative matrix $A$ as a weighted adjacency matrix of a weighted graph $G$ with vertex set $[n]$, and define a graph-theoretic quantity known as the \emph{$k$-way expansion constant}: 

\begin{definition}[Conductance and expansion]\label{defn:conductance}
Let $G=(V,w)$ be a weighted undirected graph with vertex set $V$ and non-negative edge weights $w(u,v)$. For $S\subseteq V$, the \emph{conductance} of $S$ is defined as
\begin{align*}
\phi_G(S) := \frac{w(E(S,S^c))}{w(S)} := \frac{\sum_{u\in S, v\notin S} w(u,v) }{\sum_{u\in S}\sum_{v\neq u} w(u,v)}. 
\end{align*}
For $k\ge 2$, the \emph{$k$-way expansion constant} of $G$ is defined as
\begin{align*}
\rho_G(k) := \min_{S_1,\dots,S_k\subseteq V} \max_{i\in [k]} \phi_G(S_i), 
\end{align*}
where the minimization is over all non-empty and disjoint subsets $S_1,\dots,S_k$ of $V$. 
\end{definition}

Since the channel overlap matrix $A$ is doubly stochastic, all vertices of the weighted graph $G$ have weighted degree $1$ (including self-loops). In other words, $I-A$ is a \emph{weighted Laplacian}. The following deep result, known as the \emph{higher-order Cheeger inequality} \cite{lee2014multiway}, establishes a link between the eigenvalues of the weighted Laplacian and the expansion constants of $G$. 

\begin{lemma}[Higher-order Cheeger inequality]\label{lemma:cheeger}
Let $G$ be any weighted graph, with a weighted adjacency matrix $A$ being doubly stochastic. For every $k\ge 2$, 
\begin{align*}
\frac{1-\lambda_k(A)}{2} \le \rho_G(k) \le Ck^2\sqrt{1-\lambda_k(A)}, 
\end{align*}
where $C>0$ is a universal constant. In addition, 
\begin{align*}
\rho_G(k) \le C\sqrt{(1-\lambda_{2k}(A))\log k}. 
\end{align*}
\end{lemma}

When $k=2$, the first inequality is the celebrated Cheeger inequality. Using the second inequality with a better dependence on $k$, \Cref{lemma:spectral_gap} is then a direct consequence of the inequality
\begin{align}\label{eq:combinatorial_target}
\rho_G(5k) \ge \frac{1}{4} e^{-\Diam}. 
\end{align}
The proof of \eqref{eq:combinatorial_target} is via a combinatorial argument. Fix any non-empty disjoint subsets $S_1,\dots,S_{5k}\subseteq [n]$. Since $\sum_{j\neq i} A_{ij} \le 1$ for all $i\in [n]$, it suffices to prove that
\begin{align}\label{eq:A_sum}
\max_{\ell\in [5k]} \frac{1}{|S_\ell|}\sum_{i\in S_\ell, j\notin S_\ell} A_{ij} \ge \frac{1}{4} e^{-\Diam}. 
\end{align}
The following lemma establishes a useful lower bound of $\sum_{j\notin S_\ell} A_{ij}$. 

\begin{lemma}\label{lemma:renyi-divergence}
Let $i\in S\subseteq [n]$, and $A$ be constructed from probability distributions $P_1,\dots,P_n$. For $t>0$ and $T:= \{j\notin S: D_{1/2}(P_i,P_j)\le t\}$, then
\begin{align*}
\sum_{j\notin S} A_{ij} \ge \frac{|T|}{|S|+|T|}e^{-t}. 
\end{align*}
\end{lemma}
\begin{proof}
The claim follows from the following chain of inequalities: 
  \begin{align*}
    \sum_{j \notin S} A_{ij} &= \sum_{j \notin S} \int \frac{\rmd P_i \rmd P_j}{\sum_{k=1}^n \rmd P_k} \\
    &= \int \frac{\rmd P_i \sum_{j \notin S} \rmd P_j}{\sum_{k\in S} \rmd P_k + \sum_{j \notin S} \rmd P_j } \\
    &\stepa{\ge} \int \frac{\rmd P_i \sum_{j \in T} \rmd P_j}{\sum_{k\in S} \rmd P_k + \sum_{j \in T} \rmd P_j } \\
    &\stepb{=} \frac{1}{|S| + |T|} \sum_{j\in T} \int \frac{\rmd P_i \rmd P_j}{\rmd Q} \\
    &\stepc{\ge} 
    \frac{1}{|S| + |T|} \sum_{j\in T} \pth{\int \sqrt{\rmd P_i \rmd P_j}}^2 \\
    &= \frac{1}{|S| + |T|} \sum_{j\in T} e^{-D_{1/2}(P_i,P_j)} \ge \frac{|T|}{|S| + |T|} e^{-t}, 
  \end{align*}
  where (a) follows from the increasing property of $x\mapsto \frac{ax}{b+x}$ for $a,b,x\ge 0$, (b) defines a probability measure $Q:=\frac{1}{|S|+|T|}\sum_{j\in S\cup T} P_j$, and (c) is the Cauchy--Schwarz inequality. 
\end{proof}

{\blue The next result is a combinatorial packing lemma. For $i,j\in [n]$, let $d(i,j) = D_{1/2}(P_i,P_j)$, and $N_r(i)=\{j\in [n]: d(i,j) \le r\}$ for $r\ge 0$. We call $i$ and $j$ are \emph{$r$-close} if $d(i,j)\le r$. 

\begin{lemma}\label{lemma:combinatorial}
Let $S_1,\dots,S_{5k}$ be disjoint subsets of $[n]$, ordered so that $|S_1|\le \dots\le |S_{5k}|$. Suppose $[n]$ is partitioned into $k$ cells $\calC_1,\dots,\calC_k$, and every two indices in the same cell are $r$-close. Then for some $\ell \in [5k]$, 
\begin{align*}
|\{i\in S_\ell: |N_r(i)|\ge 2|S_\ell|\}| \ge \frac{|S_\ell|}{2}. 
\end{align*}
\end{lemma}
\begin{proof}
Suppose not. Let $B_\ell \subseteq S_\ell$ be the bad indices:
\begin{align*}
B_{\ell} := \{i\in S_\ell: |N_r(i)|< 2|S_\ell|\}, 
\end{align*}
so $|B_\ell| > |S_\ell|/2$ for every $\ell\in [5k]$. For each $\calC_m$ meeting some $B_\ell$, let
\begin{align*}
a_m := \min\{\ell: \calC_m \cap B_{\ell} \neq \varnothing \}. 
\end{align*}
If $i\in \calC_m \cap B_{a_m}$, then $\calC_m\subseteq N_r(i)$, hence
\begin{align}\label{eq:property}
|\calC_m| \le |N_r(i)| < 2|S_{a_m}|. 
\end{align}
Now let $M_t:=\{m: a_m\le 4t+1\}$. The central claim is that $|M_t|\ge t+1$ for every $t=0,1,\dots,k$. We prove it by induction on $t$. The base case $t=0$ is obvious. Now suppose that $|M_s|\ge s+1$ for all $s=0,\dots,t-1$, but $|M_t|\le t$. Every bad point in $B_1\cup \dots\cup B_{4t+1}$ lies in $\cup_{m\in M_t} \calC_m$, so by \eqref{eq:property}, 
\begin{align*}
\sum_{\ell=1}^{4t+1} |B_\ell| \le \sum_{m\in M_t} |\calC_m| \le 2\sum_{m\in M_t} |S_{a_m}|. 
\end{align*}
Since $|M_s|\ge s+1$ for all $s=0,\dots,t-1$ and $|M_t|\le t$, the right-hand side is at most
\begin{align*}
2\sum_{j=1}^t |S_{4j-3}| \le \frac{1}{2}\sum_{\ell=1}^{4t+1} |S_{\ell}|. 
\end{align*}
But the bad-set assumption gives $\sum_{\ell=1}^{4t+1} |B_\ell| > \frac{1}{2}\sum_{\ell=1}^{4t+1} |S_{\ell}|$, a contradiction. Therefore, the claim $|M_t|\ge t+1$ holds for every $t=0,1,\dots,k$; taking $t=k$ yields $|M_k|\ge k+1$, a contradiction because there are only $k$ cells. 
\end{proof}

Next we prove \eqref{eq:A_sum}. By the definition of R\'enyi partition diameter $\Diam$, the condition of \Cref{lemma:combinatorial} holds with $r=\Diam$. By \Cref{lemma:combinatorial}, we can take $\ell\in [5k]$ such that 
\begin{align*}
|\{i\in S_\ell: |N_r(i)|\ge 2|S_\ell|\}| \ge \frac{|S_\ell|}{2}. 
\end{align*}
For every $i$ in the above set, applying \Cref{lemma:renyi-divergence} with $t=r$ gives 
\begin{align*}
\sum_{j\notin S_{\ell}} A_{ij} \ge \frac{|N_r(i)\backslash S_{\ell}|}{|S_{\ell}|+|N_r(i)\backslash S_{\ell}|}e^{-r} \ge \frac{1}{2}e^{-r},
\end{align*}
for $|N_r(i)\backslash S_{\ell}|\ge |N_r(i)| - |S_{\ell}| \ge |S_{\ell}|$. Finally,
\begin{align*}
\frac{1}{|S_{\ell}|}\sum_{i\in S_{\ell}}\sum_{j\notin S_{\ell}} A_{ij} \ge \frac{1}{|S_{\ell}|}\cdot \frac{|S_{\ell}|}{2}\cdot \frac{1}{2}e^{-r} = \frac{1}{4}e^{-r}, 
\end{align*}
which is the target inequality \eqref{eq:A_sum}. 
}

\subsection{Proof of the dimension-dependent upper bound}
{\blue This section proves \Cref{thm:main_upper_dimdep} via a refined upper bound of \eqref{eq:A-upper}, shown in \cite[Lemma 5.3]{han2024approximate}:
\begin{align*}
\chi^2(\bP_n \| \bQ_n) + 1  \le &\sum_{\ell=0}^n \sum_{\substack{\ell_2+\dots+\ell_n=\ell \\ \ell_2,\dots,\ell_n\ge 0}} \lambda_2(A)^{\ell_2}\cdots \lambda_n(A)^{\ell_n} \le \prod_{i=2}^n \sum_{\ell_i=0}^n \lambda_i(A)^{\ell_i} \\\le &\prod_{i=2}^n \min\sth{\frac{1}{1-\lambda_i(A)},n+1}, 
\end{align*}
where $A\in \bR^{n\times n}$ is the channel overlap matrix. As a result, 
\begin{align}\label{eq:A-upper-refined}
\log(\chi^2(\bP_n \| \bQ_n) + 1)&\le \sum_{i=2}^n \min\sth{ \log \frac{1}{1-\lambda_i(A)}, \log (n+1) }  =: \sum_{i=2}^n \psi_n(\lambda_i(A)). 
\end{align}
Next we apply \Cref{lemma:spectral_gap} to the truncated upper bound \eqref{eq:A-upper-refined}. 

We first control the eigenvalues \(\lambda_2(A),\ldots,\lambda_9(A)\).
By~\cite[Lemma 5.2]{han2024approximate}, for every \(i\ge2\), \(
1-\lambda_i(A)\ge e^{-D_1(\mathcal P)}
\). 
Consequently,
\begin{equation}
\sum_{i=2}^{9}\psi_n(\lambda_i(A))
\le
8\min\{D_1(\mathcal P),\log(n+1)\}. 
\label{eq:first-eight-eigenvalues}
\end{equation}
It remains to control the eigenvalues with indices $i\ge 10$.
For \(k\ge1\), define $g_k := \max_{\ell\le k} ce^{-2\sfD_{\ell}(\calP)}/\log(5\ell)$, 
where \(c>0\) is the universal constant in Lemma~\ref{lemma:spectral_gap}.
Decreasing \(c\) if necessary, we may assume that \(g_k\le 1/2\)
for every \(k\). Since \(g_k\) is nondecreasing, Lemma~\ref{lemma:spectral_gap}
and the monotonicity of the eigenvalues imply
\begin{equation}
0\le
\lambda_{10k}(A),\ldots,\lambda_{10k+9}(A)
\le 1-g_k,
\qquad k\ge1.
\label{eq:block-eigenvalue-cap-finite}
\end{equation}
Moreover, as in \eqref{eq:constraint_2},
\begin{equation}
\sum_{i=10}^n\lambda_i(A)
\le \trace(A)-1
\le \Capa.
\label{eq:trace-budget-finite}
\end{equation}

We next optimize $\sum_{i=10}^n \psi_n(\lambda_i(A))$ subject to \eqref{eq:block-eigenvalue-cap-finite} and \eqref{eq:trace-budget-finite}. First, if $\lambda_i(A)>\frac{n}{n+1}$, we may replace $\lambda_i(A)$ by $\min\{\lambda_i(A),\frac{n}{n+1}\}$ without decreasing the objective \eqref{eq:A-upper-refined}, while keeping \eqref{eq:block-eigenvalue-cap-finite} and \eqref{eq:trace-budget-finite} both valid. Therefore, we may restrict attention to variables in \([0,\frac{n}{n+1}]\), on
which \(x\mapsto\psi_n(x)=\log(1/(1-x))\) is increasing and convex. Then by the same analysis in the proof of \Cref{thm:main_upper}, the maximizer has the following structure: for some block index
\(k_0\), 
\begin{align*}
x_{10k}^\star, \dots, x_{10k+9}^\star \begin{cases}
   = 1-\max\{g_k, \frac{1}{n+1}\} & \text{for } k < k_0 \\
   \in [0,1-\max\{g_k, \frac{1}{n+1}\}] & \text{for } k = k_0 \\
   = 0 & \text{for } k > k_0
\end{cases}, \qquad \text{with }\sum_{i=10}^n x_i^\star = \Capa. 
\end{align*}
The cutoff \(k_0\) is controlled by the trace budget. Since $\max\{g_k,\frac{1}{n+1}\}\le \frac{1}{2}$, 
\[
\Capa
\ge
\sum_{i=10}^{10k_0-1} x_i^\star
\ge
5(k_0-1),
\]
so that $k_0\le 1+\Capa$. Next, since
\begin{align*}
\min\left\{\log\frac1{g_k},\log(n+1)\right\}
 \le
&2\min\{D_k(\mathcal P),\log(n+1)\}
+
\log\frac{\log(5k)}{c} \\
\lesssim &
\min\{D_k(\mathcal P),\log n\}
+
\log_+\log(5k), 
\end{align*}
we obtain
\begin{align*}
\sum_{i=10}^n\psi_n(x_i^\star) 
\lesssim&
\sum_{k=1}^{\lfloor\Capa\rfloor+1}
\min\{D_k(\mathcal P),\log n\}
+
\sum_{k=1}^{\lfloor\Capa\rfloor+1}
\log_+\log(5k)
\nonumber\\
\lesssim&
\sum_{k=1}^{\lfloor\Capa\rfloor+1}
\min\{D_k(\mathcal P),\log n\}
+ 
(\Capa+1)\log^+\log \Capa.
\end{align*}
Together with \eqref{eq:first-eight-eigenvalues}, the claimed bound holds for $\log(\chi^2(\bP_n\|\bQ_n)+1)$. Replacing $n$ by any $m\le n$, and taking the supremum over $P_1,\dots,P_m\in \calP$, \Cref{thm:main_upper_dimdep} follows. 
}

\subsection{Lower bounds via saddle point analysis}
The key ingredient in the proof of \Cref{thm:main_lower} is the following characterization of a limiting matrix permanent.

\begin{lemma}\label{lemma:permanent}
Let $A\in \bR^{n\times n}$ be a symmetric doubly stochastic matrix with strictly positive entries, and $J_m$ be the $m\times m$ all-ones matrix. Then
\begin{align*}
\liminf_{m\to\infty} \frac{(mn)^{mn}}{(mn)!} \Perm\pth{A\otimes \frac{J_m}{m}} \ge \prod_{k=2}^n \frac{1}{\sqrt{1-\lambda_k(A)^2}}. 
\end{align*}
\end{lemma}
\begin{remark}
We conjecture that this is also the limit as $m \to \infty$, although establishing it would require a more delicate saddle-point analysis. For the proof of \Cref{thm:main_lower}, however, we only require the lower bound direction.
\end{remark}

We show that \Cref{lemma:permanent} implies \Cref{thm:main_lower}. First, if we can find mutually singular distributions $P_1, P_2\in \calP$, then \cite[Lemma 5.1]{han2025best} shows that $\chi^2(\calP) = \infty$, and there is nothing to prove. Therefore, for given $P_1,\dots,P_n\in \calP$, we may assume that the channel overlap matrix $A$ has strictly positive entries. Now let each $P_i$ appear $m$ times in the permutation mixture $\bP_{mn}$, so that the resulting channel overlap matrix $A_m\in \bR^{mn\times mn}$ becomes $A_m = A\otimes \frac{J_m}{m}$. By \eqref{eq:permanent}, 
\begin{align*}
1 + \chi^2(\bP_{mn} \| \bQ_{mn}) = \frac{(mn)^{mn}}{(mn)!} \Perm\pth{A\otimes \frac{J_m}{m}}. 
\end{align*}
Letting $m\to\infty$, \Cref{lemma:permanent} gives the first lower bound of \Cref{thm:main_lower}. For the second lower bound, since $A$ is doubly stochastic, $I-A^2$ is a matrix with zero row and column sums. Then by Kirchhoff's matrix-tree theorem \cite{kirchhoff1847ueber}, 
\begin{align*}
  \prod_{k=2}^n (1-\lambda_k(A)^2) = \prod_{k=1}^{n-1} \lambda_k(I-A^2) = {\blue n}\det((I-A^2)_{(1,1)})
\end{align*}
where $M_{(1,1)}$ denotes the matrix obtained by removing the first row and column of $M$. Next, since $I-A^2$ is PSD by Perron--Frobenius, so is the principal submatrix $(I-A^2)_{(1,1)}$. Finally, Hadamard's inequality yields
\begin{align*}
\det((I-A^2)_{(1,1)}) \le \prod_{k=2}^n (I-A^2)_{kk} = \prod_{k=2}^n (1-(A^2)_{kk}), 
\end{align*}
completing the proof of the second inequality. 

\begin{proof}[Proof of \Cref{lemma:permanent}]
{\blue We first outline the proof before turning to the details. We express the permanent $\Perm\pth{A\otimes \frac{J_m}{m}}$ as a sum over $n\times n$ contingency tables whose row and column sums are all equal to $m$. As $m\to\infty$, the sum over normalized contingency tables converges to an integral with respect to the uniform measure on the space of doubly stochastic matrices. We then apply Laplace's approximation to obtain a lower bound on the limiting integral. The resulting expression involves the determinant of the Hessian matrix appearing in Laplace's approximation, which we compute exactly through linear-algebraic manipulations.
}

In the first step, we establish the following identity: 
\begin{align}\label{eq:contingency}
&\Perm\pth{A\otimes \frac{J_m}{m}} = (m!)^{n}  \sum_{X\in \calT_m} \prod_{i=1}^n \binom{m}{X_{i1}, \dots, X_{in}}\prod_{j=1}^n \left(\frac{A_{ij}}{m}\right)^{X_{ij}}, 
\end{align}
where $\calT_m$ denotes the class of $n\times n$ contingency tables with all rows and columns summing to $m$: 
\begin{align*}
\calT_m = \bsth{ X\in \naturals^{n\times n}: \sum_{i=1}^n X_{ij} = m, \forall j\in [n]; \sum_{j=1}^n X_{ij} = m, \forall i\in [n]}. 
\end{align*}
To show \eqref{eq:contingency}, note that $A\otimes \frac{J_m}{m}$ consists of $n^2$ blocks of $m\times m$ matrices, where the $(i,j)$-th block has the same entry $A_{ij}/m$. For a permutation $\pi\in S_{mn}$, let $X_{ij} =|\pi(B_i) \cap B_j| $ $:= |\pi(\{(i-1)m+1,\dots,im\})\cap \{(j-1)m+1,\dots,jm\}|$ be the number of rows in block $i$ that are matched to columns within block $j$. Clearly $X\in \calT_m$. Next we count the number of permutations $\pi\in S_{mn}$ that map to a given contingency table $X$. To this end, we first enumerate over the membership of all rows, i.e. the indices of column blocks each row is mapped into. The number of memberships giving rise to a given contingency table $X$ is
\begin{align*}
    \prod_{i=1}^n \binom{m}{X_{i1}, \dots, X_{in}}. 
\end{align*}
Second, for a given membership, the number of ways to assign column indices within blocks is $(m!)^n$, by enumerating over column blocks. This establishes \eqref{eq:contingency}. 

Next we apply Stirling's approximation $m!=(1+o_m(1))\sqrt{2\pi m}(\frac{m}{e})^m$ for large $m$, where we will frequently use $o_m(1)$ to denote terms that go to $0$ as $m\to\infty$. Starting from \eqref{eq:contingency}, 
\begin{align*}
    \frac{(mn)^{mn}}{(mn)!}\text{Perm}(A\otimes \frac{J_m}{m}) 
      =& \frac{(mn)^{mn}}{(mn)!} \cdot (m!)^{n}  \sum_{X\in \calT_m} \prod_{i=1}^n \binom{m}{X_{i1}, \dots, X_{in}}\prod_{j=1}^n \left(\frac{A_{ij}}{m}\right)^{X_{ij}} \\
      =& \frac{(mn)^{mn}}{(mn)!} \cdot (m!)^{2n}  \sum_{X\in \calT_m} \prod_{i,j=1}^n \frac{1}{X_{ij}!}\left(\frac{A_{ij}}{m}\right)^{X_{ij}}\\
      \stepa{=}& \frac{n^{mn}}{(mn)!} \cdot (m!)^{2n}  \sum_{X\in \calT_m} \prod_{i,j=1}^n \frac{A_{ij}^{X_{ij}}}{X_{ij}!} \\
      =& \frac{n^{mn}}{\binom{mn}{m, \dots, m}} \cdot (m!)^{n} \sum_{X\in \calT_m} \exp\bqth{ \sum_{i,j=1}^n \bpth{X_{ij} \log A_{ij} - \log (X_{ij}!) }} \\
      \stepb{=}& (1+o_m(1))\frac{(\sqrt{2\pi m})^n}{\sqrt{2\pi mn}} \cdot (m!)^n \sum_{X\in \calT_m} \exp\bqth{ \sum_{i,j=1}^n \bpth{X_{ij} \log A_{ij} - \log (X_{ij}!) }},
    \end{align*}
where (a) uses the definition of contingency table that $\sum_{i,j} X_{ij}=mn$, and (b) follows from Stirling's approximation (note that $n$ is a fixed number, so $(1+o_m(1))^n = 1+o_m(1)$). Next we lower bound the above quantity by summing over a subset of $\calT_m$. Since $A$ has strictly positive entries, we have $A_{\min} := \min_{i,j} A_{ij} > 0$. Choosing any sequence $\varepsilon_m\in (0, A_{\min}/2)$ with $\sqrt{m}\varepsilon_m\to B$ (where $B$ is any fixed constant), define the subset
\begin{align*}
\calT_m' := \bsth{X\in \calT_m: |X_{ij} - mA_{ij}|\le m\varepsilon_m, \forall i,j\in [n]} \subseteq \calT_m. 
\end{align*}
By the choice of $\varepsilon_m<\frac{A_{\min}}{2}$, for all $X\in \calT_m'$ we have $X_{ij}\ge \frac{mA_{\min}}{2} = \Omega(m)$. Therefore, we may apply Stirling's approximation to $\log(X_{ij}!)=X_{ij}(\log X_{ij} - 1) + \frac{1}{2}\log(2\pi X_{ij})+o_m(1)$ to obtain
\begin{align*}
    &\frac{(mn)^{mn}}{(mn)!}\text{Perm}(A\otimes \frac{J_m}{m}) \\
      &\ge (1+o_m(1))\frac{(\sqrt{2\pi m})^n}{\sqrt{2\pi mn}} \cdot (m!)^n \sum_{X\in \calT_m'} \exp\bqth{ \sum_{i,j=1}^n \bpth{X_{ij} \log A_{ij} - \log (X_{ij}!) }} \\
      &= (1+o_m(1))\frac{(\sqrt{2\pi m})^n}{\sqrt{2\pi mn}} \cdot (m!)^n  \sum_{X\in \calT_m'} \exp\bqth{ \sum_{i,j=1}^n \bpth{ X_{ij}\log \frac{A_{ij}}{X_{ij}} +X_{ij}-\frac{1}{2}\log(2\pi X_{ij}) }} \\
      &\stepc{=} (1+o_m(1))\frac{(\sqrt{2\pi m})^n}{\sqrt{2\pi mn}} \cdot  \frac{(m!)^n e^{mn}}{\prod_{i,j=1}^n \sqrt{2\pi m A_{ij}}} \sum_{X\in \calT_m'} \exp\bqth{ \sum_{i,j=1}^n X_{ij}\log \frac{A_{ij}}{X_{ij}}  } \\
      &\stepd{=} (1+o_m(1))\frac{1}{\sqrt{n\prod_{i,j=1}^n A_{ij}}} \cdot  \frac{1}{(\sqrt{2\pi m})^{(n-1)^2}} \sum_{X\in \calT_m'} \exp\bqth{ \sum_{i,j=1}^n X_{ij}\log \frac{mA_{ij}}{X_{ij}}  }, 
    \end{align*}
    where (c) uses that $\log(X_{ij})=\log(mA_{ij})+o_m(1)$ by definition of $\calT_m'$ and $\varepsilon_m=o_m(1)$, and (d) follows from Stirling's approximation to $m!$, the identity $\sum_{i,j} X_{ij}=mn$, and rearranging. 

    Next, using saddle point analysis, the following lemma establishes an asymptotic lower bound of the remaining sum. 
    \begin{lemma}\label{lemma:saddle-point}
    Let $\Omega$ be the set of all positive doubly stochastic matrices, and for $z\in \Omega$, define
    \begin{align*}
    f(z) = \sum_{i,j=1}^n z_{ij}\log \frac{z_{ij}}{A_{ij}}. 
    \end{align*}
    Then
    \begin{align*}
    \liminf_{B\to\infty}\liminf_{m\to\infty} \frac{1}{(\sqrt{2\pi m})^{(n-1)^2}} \sum_{X\in \calT_m'} \exp\bqth{ \sum_{i,j=1}^n X_{ij}\log \frac{mA_{ij}}{X_{ij}}  }  \ge \frac{1}{\sqrt{\det(\nabla^2 f(A))}}, 
    \end{align*}
    where $\nabla^2 f$ is the constrained Hessian of $f$ by treating $\{z_{ij}\}_{i,j=1}^{n-1}$ as the only free variables. 
    \end{lemma}
    \begin{proof}
    Let $X_{ij} = mz_{ij}$, then 
\begin{align*}
    \sum_{X\in \calT_m'} \exp\bqth{ \sum_{i,j=1}^n X_{ij}\log \frac{mA_{ij}}{X_{ij}}  } = \sum_{z\in \frac{1}{m}\calT_m'} \exp\pth{ -mf(\widetilde{z}) }. 
\end{align*}
To avoid ambiguity, now we treat $f$ as a function of the $(n-1)^2$ free variables $\{z_{ij}\}_{i,j=1}^{n-1}$; we write $\widetilde{z}$ and $\widetilde{A}$ to denote the restricted matrices. By direct calculation, $f(\widetilde{A}) = 0$, $\nabla f(\widetilde{A})=0$, and $\nabla^3 f(\widetilde{z})$ is an $(n-1)^2\times (n-1)^2 \times (n-1)^2$ tensor with a bounded norm whenever $\min_{i,j} z_{ij} = \Omega(1)$. Expanding $f$ in a third-order Taylor series of $\widetilde{A}$, we obtain
\begin{align*}
f(\widetilde{z}) &= \frac{1}{2}(\widetilde{z}-\widetilde{A})^\top \nabla^2 f(\widetilde{A}) (\widetilde{z}-\widetilde{A}) + O(\|\widetilde{z}-\widetilde{A}\|^3) \\ &=  \frac{1}{2}(\widetilde{z}-\widetilde{A})^\top \nabla^2 f(\widetilde{A}) (\widetilde{z}-\widetilde{A}) + O(\varepsilon_m^3)
\end{align*}
uniformly for $z\in \frac{1}{m}\calT_m'$, where the properties of $\calT_m'$ guarantee that $\min_{i,j} z_{ij} = \Omega(1)$ and $\|z-A\|=O(\varepsilon_m)$. Since $m\varepsilon_m^3\to 0$, it follows that
\begin{align*}
&\frac{1}{(\sqrt{2\pi m})^{(n-1)^2}} \sum_{X\in \calT_m'} \exp\bqth{ \sum_{i,j=1}^n X_{ij}\log \frac{mA_{ij}}{X_{ij}}  } \\&\ge \frac{1+o_m(1)}{(\sqrt{2\pi m})^{(n-1)^2}} \sum_{z\in \frac{1}{m}\calT_m'} \exp\pth{-\frac{m}{2}(\widetilde{z}-\widetilde{A})^\top \nabla^2 f(\widetilde{A}) (\widetilde{z}-\widetilde{A})}. 
\end{align*}
Make the change of variables $\xi = \sqrt{m}(z-A)$, and let $\Omega_m$ be the set of all possible values of $\widetilde{\xi}$ as $z$ ranges over $\frac{1}{m}\calT_m'$. This set is contained in a lattice of determinant $(\frac{1}{\sqrt{m}})^{(n-1)^2}$. Since $\sqrt{m}\varepsilon_m\to B$, the lattice polytope (i.e. the convex hull of $\Omega_m$) converges to the set $S(B)$, which is the projection onto the first $(n-1)^2$ coordinates of
\begin{align*}
\bsth{ \xi\in \bR^{n\times n}: \|\xi\|_\infty \le B, \sum_{i=1}^n \xi_{ij}=0, \forall j\in [n]; \sum_{j=1}^n \xi_{ij}=0, \forall i\in [n] }. 
\end{align*}
Therefore
\begin{align*}
&\frac{1}{(\sqrt{2\pi m})^{(n-1)^2}}\sum_{z\in \frac{1}{m}\calT_m'} \exp\pth{-\frac{m}{2}(\widetilde{z}-\widetilde{A})^\top \nabla^2 f(\widetilde{A}) (\widetilde{z}-\widetilde{A})} \\ &= \frac{1}{(\sqrt{2\pi m})^{(n-1)^2}}\sum_{\widetilde{\xi}\in \Omega_m} \exp\pth{-\frac{1}{2}\widetilde{\xi}^\top \nabla^2 f(\widetilde{A}) \widetilde{\xi}}
\end{align*}
is a finite Riemann sum on the set $S(B)$, which converges to 
\begin{align*}
\int_{S(B)} \frac{1}{(\sqrt{2\pi})^{(n-1)^2}} \exp\pth{-\frac{1}{2}\widetilde{\xi}^\top \nabla^2 f(\widetilde{A}) \widetilde{\xi}} \rmd \widetilde{\xi}. 
\end{align*}
Finally, since $S(B)\to \mathbb{R}^{(n-1)^2}$ as $B\to\infty$, the above Gaussian integral converges to $(\det(\nabla^2 f(\widetilde{A})))^{-1/2}$, as desired. 
\end{proof} 
   
Since the choice of the constant $B>0$ in the definition of $\varepsilon_m$ is arbitrary, \Cref{lemma:saddle-point} implies \Cref{lemma:permanent} provided that the constrained Hessian $H=\nabla^2 f(A)$ has determinant
\begin{align}\label{eq:constrained_Hessian}
    \det(H) = \frac{1}{n\prod_{i,j=1}^n A_{ij}}\prod_{k=2}^n (1-\lambda_k(A)^2). 
\end{align}
The proof of \eqref{eq:constrained_Hessian} requires careful linear algebraic manipulations. Defining all other variables
    \begin{align*}
    \begin{cases}
    z_{in} = 1-\sum_{j<n}z_{ij}, 
    \\
    z_{nj} = 1-\sum_{i<n}z_{ij}, 
    \\
    z_{nn} = 1-\sum_{i<n}\sum_{j<n}z_{ij},
    \end{cases}
    \end{align*}
as functions of the free variables $\{z_{ij}\}_{i,j=1}^{n-1}$, we can compute the Hessian as
    \begin{align*}
    H_{ij,i'j'}=\frac{\partial^2 f(A)}{\partial z_{ij}\,\partial z_{i'j'}}=
    \frac{\delta_{ii'}\,\delta_{jj'}}{A_{ij}}
    +\frac{\delta_{ii'}}{A_{in}}
    +\frac{\delta_{jj'}}{A_{nj}}
    +\frac{1}{A_{nn}},\quad i,j,i',j'\in [n-1]. 
    \end{align*}
    Alternatively, we can write
    \begin{align*}
    H =
    D_{0}
   +
    \sum_{i=1}^{n-1}\frac1{A_{in}}\,r_{i}r_{i}^{\top}
    +
    \sum_{j=1}^{n-1}\frac1{A_{nj}}\,c_{j}c_{j}^{\top}
    +
    \frac1{A_{nn}}\,\mathbf1_{(n-1)^2}\,\mathbf1_{(n-1)^2}^{\top},
    \end{align*}
    where
    \begin{align*}
      D_{0}  = \operatorname{diag}(\frac{1}{A_{ij}}), 1 \leq i,j \leq n-1, \quad r_i = e_i \otimes \mathbf1_{n-1}, \quad c_j = \mathbf1_{n-1} \otimes e_j.
    \end{align*}
   Here $e_i \in \bR^{n-1}$ is the basis vector in $\bR^{n-1}$. A more compact expression is
    \begin{align*}
    H = D_{0} + U D_1 U^\top
    \end{align*}
    where
    \begin{align*}
    U = \left[r_1, \dots r_{n-1}, c_1, \dots c_{n-1}, \mathbf1_{(n-1)^2}\right], \quad D_1 = \operatorname{diag}(\frac{1}{A_{1n}}, \dots, \frac{1}{A_{n-1,n}}, \frac{1}{A_{n1}}, \dots \frac{1}{A_{n,n-1}}, \frac{1}{A_{nn}}).
    \end{align*}
    
    Using the above expression, the matrix determinant lemma yields
    \begin{align*}
      \det(H) 
      &= \det\pth{D_{0} + U D_1 U^\top } \\& = \det(D_0) \det(D_1) \det\pth{D_1^{-1} + U^\top D_{0}^{-1} U }.
    \end{align*}
        Next we explore the structure of the matrix $ N := D_1^{-1} + U^\top D_{0}^{-1} U $: 
\begin{align*}
N
&= D_1^{-1}
+
\begin{bmatrix}
\mathbf1^\top & \mathbf0^\top & \dots & \mathbf0^\top \\
\vdots & \vdots & \ddots & \vdots \\
\mathbf0^\top & \mathbf0^\top & \dots & \mathbf1^\top \\ 
e_1^\top& e_1^\top & \dots & e_1^\top  \\
\vdots & \vdots & \ddots & \vdots \\
e_{n-1}^\top & e_{n-1}^\top &\dots & e_{n-1}^\top \\ 
\mathbf1^\top & \mathbf1^\top & \dots & \mathbf1^\top
\end{bmatrix}
D_0^{-1}
\begin{bmatrix}
\mathbf1 & \mathbf0 & \cdots & \mathbf0
  & e_1 & e_2 & \cdots & e_{n-1} & \mathbf1\\
\mathbf0 & \mathbf1 & \cdots & \mathbf0
  & e_1 & e_2 & \cdots & e_{n-1} & \mathbf1\\
\vdots & \vdots & \ddots & \vdots
  & \vdots & \vdots & \ddots & \vdots & \vdots\\
\mathbf0 & \mathbf0 & \cdots & \mathbf1
  & e_1 & e_2 & \cdots & e_{n-1} & \mathbf1
\end{bmatrix} \\
&= D_1^{-1}
+
\begin{bmatrix}
\begin{matrix}
s_1 & & \\
& \ddots & \\
& & s_{n-1}
\end{matrix}
&
\widetilde A
&
\begin{matrix}
s_1\\
\vdots\\
s_{n-1}
\end{matrix}
\\
\widetilde A
&
\begin{matrix}
s_1 & & \\
& \ddots & \\
& & s_{n-1}
\end{matrix}
&
\begin{matrix}
s_1\\
\vdots\\
s_{n-1}
\end{matrix}
\\
\begin{matrix}
s_1 & \cdots & s_{n-1}
\end{matrix}
&
\begin{matrix}
s_1 & \cdots & s_{n-1}
\end{matrix}
&
\displaystyle\sum_{i=1}^{n-1}s_i
\end{bmatrix}
=
\begin{bmatrix}
I_{n-1}
&
\widetilde A
&
\begin{matrix}
s_1\\
\vdots\\
s_{n-1}
\end{matrix}
\\
\widetilde A
&
I_{n-1}
&
\begin{matrix}
s_1\\
\vdots\\
s_{n-1}
\end{matrix}
\\
\begin{matrix}
s_1 & \cdots & s_{n-1}
\end{matrix}
&
\begin{matrix}
s_1 & \cdots & s_{n-1}
\end{matrix}
&
A_{nn}+\displaystyle\sum_{i=1}^{n-1}s_i
\end{bmatrix},
\end{align*}
where $\widetilde{A}=(A_{ij})_{i,j=1}^{n-1}$ is the upper left corner of $A$, and $s_i = \sum_{j=1}^{n-1} A_{ij}$ is the $i$-th row (or column) sum of $\widetilde{A}$. To proceed, we apply the following elementary row and column operations to $N$: 
    \begin{enumerate}
      \item First, we subtract the last row of $N$ by the sum of $(n,n+1,\dots,2n-2)$-th rows to obtain
      \begin{align*}
          N_1 := \begin{bmatrix}
        I_{n-1} & \widetilde{A} & \begin{matrix}
            s_1 \\
            \vdots \\
            s_{n-1}
        \end{matrix} \\
        \widetilde{A} & I_{n-1} & \begin{matrix}
            s_1 \\
            \vdots \\
            s_{n-1}
        \end{matrix} \\
        \begin{matrix}
            0 & \cdots & 0
        \end{matrix} & \begin{matrix}
            -A_{n1} & \cdots & -A_{n,n-1}
        \end{matrix} & A_{nn}
        \end{bmatrix}; 
      \end{align*}
      \item Next, we subtract the last column of $N_1$ by the sum of $(n,n+1,\dots,2n-2)$-th columns to obtain
      \begin{align*}
        N_2 := \begin{bmatrix}
        I_{n-1} & \widetilde{A} & \begin{matrix}
            0 \\
            \vdots \\
            0
        \end{matrix} \\
        \widetilde{A} & I_{n-1} & \begin{matrix}
            -A_{n1} \\
            \vdots \\
            -A_{n,n-1}
        \end{matrix} \\
        \begin{matrix}
            0 & \cdots & 0
        \end{matrix} & \begin{matrix}
            -A_{n1} & \cdots & -A_{n,n-1}
        \end{matrix} & 1
        \end{bmatrix};
      \end{align*}
      \item Finally, we multiply $-1$ to the first $(n-1)$ rows and columns of $N_2$, and move the last row and column to the middle to obtain
      \begin{align*}
          N_3 := \begin{bmatrix}
              I_n & -A \\
              -A & I_n
          \end{bmatrix} \text{ with last row and column removed}. 
      \end{align*}
    \end{enumerate}
    Since the above elementary operations do not change the determinant of $N$, we conclude that $\det(N)$ is the $(2n,2n)$-th cofactor of the matrix 
    \begin{align*}
        M = \begin{bmatrix}
              I_n & -A \\
              -A & I_n
          \end{bmatrix}. 
    \end{align*}
 Since $M$ has zero row and column sums, Kirchhoff matrix-tree theorem \cite{kirchhoff1847ueber} yields
    \begin{align*}
    \text{all cofactors of }M = \frac{\prod_{k=1}^{2n-1} \lambda_k(M)}{2n},
    \end{align*}
    where $\lambda_1(M),\dots,\lambda_{2n-1}(M),0$ are all eigenvalues of $M$. On the other hand, since $A$ is symmetric, all eigenvalues of $M$ are $\{1\pm \lambda_k(A)\}_{k=1}^n$. Using $\lambda_1(A) = 1$, we conclude that
    \begin{align*}
    \det(N) &= \frac{2\prod_{k=2}^n (1+\lambda_k(A))(1-\lambda_k(A))}{2n} = \frac{\prod_{k=2}^n (1-\lambda_k(A)^2)}{n}, 
    \end{align*}
    and
    \begin{align*}
    \det(H) &= \det(D_0)\det(D_1)\det(N) = \frac{1}{n\prod_{i,j=1}^n A_{ij}}\prod_{k=2}^n (1-\lambda_k(A)^2). 
    \end{align*}
    This is the desired identity \eqref{eq:constrained_Hessian}. 
\end{proof}

%% file: Application_To_Compound_Decision_Problems.tex
\section{Applications to compound decision problems}\label{sec:EB}
\subsection{Preliminaries}\label{subsec:EB_preliminary}
In this section, we review the Bayesian interpretations and some alternative definitions of separable and PI oracles in the compound decision literature.

To understand the separable oracle $\thetahatS$ in \eqref{eq:separable_oracle}, consider the following Bayesian experiment: first, draw $\theta \sim G_n := \frac{1}{n} \sum_{i=1}^n \delta_{\theta_i}$, the empirical distribution of $\theta_1, \dots, \theta_n$, and then draw $X \sim P_\theta$. The current data-generating scheme $X_i \sim P_{\theta_i}$ is related to this Bayesian experiment via the \emph{fundamental theorem of compound estimation} \cite{copas1969compound,zhang2003compound}, which states that
\begin{align}\label{eq:fundamental-theorem}
\bE\bqth{\frac{1}{n}\sum_{i=1}^n f(\theta_i, X_i)} = \bE_{G_n}\qth{f(\theta,X)}
\end{align}
for any measurable function $f$. When $f(\theta, X) = (\theta - \widehat{\theta}(X))^2$ with a separable estimator $\widehat{\theta}$, the LHS of \eqref{eq:fundamental-theorem} equals the normalized MSE of $\widehat{\theta}$. Thus, the separable oracle $\thetahatS$ is the Bayes optimal estimator under the above Bayesian experiment, i.e. $\thetahatS(X) = \bE_{G_n}[\theta|X]$, or more explicitly,
\begin{align}\label{eq:separable_oracle_explicit}
\thetahatS_i = \frac{\sum_{j=1}^n \theta_j f_{\theta_j}(X_i)}{\sum_{j=1}^n f_{\theta_j}(X_i)}, \qquad i\in [n],
\end{align}
where $f_\theta$ is the density of $P_\theta$ with respect to some dominating measure.

To understand the PI oracle $\thetahatPI$ in \eqref{eq:PI_oracle}, we use an equivalent formulation of the compound regret in \eqref{eq:regret_PI}. By \cite[Appendix A]{han2025best}, 
\begin{align}\label{eq:PI_regret_equiv}
\regPI(\widehat{\theta}) = \sup_\theta \qth{ \MSE(\theta, \widehat{\theta}) - \min_{\widehat{\theta}^\star}\max_{\pi\in S_n} \MSE(\pi(\theta), \widehat{\theta}^\star) },
\end{align}
where the inner minimization is taken over all estimators (oracle) $\widehat{\theta}^\star$ which may depend on the ground truth $\theta$. In other words, this is the optimal estimator that \emph{knows $\theta$ up to permutation}. The minimizer can be characterized via a \emph{postulated Bayes model} in which $\pi \sim \mathrm{Unif}(S_n)$, $\widetilde{\theta}_i = \theta_{\pi(i)}$, and $X_i \sim P_{\widetilde{\theta}_i}$ given $\widetilde{\theta}$ \cite{weinstein2021permutation}. Under the quadratic loss, the Bayes estimator in the postulated Bayes model is
\begin{align*}
\widehat{\theta}^\star_i = \bE[\widetilde{\theta}_i | X^n], \qquad i\in [n]. 
\end{align*}
This estimator is clearly permutation-invariant and achieves the same MSE for all $\pi \in S_n$, making it minimax and thus the inner minimizer in \eqref{eq:PI_regret_equiv}. Returning to the original formulation in \eqref{eq:regret_PI}, we obtain $\thetahatPI_i = \bE[\widetilde{\theta}_i|X^n]$, or more explicitly,
\begin{align}\label{eq:PI_oracle_explicit}
\thetahatPI_i = \frac{\sum_{\pi\in S_n} \theta_{\pi(i)} \prod_{j=1}^n f_{\theta_{\pi(j)}}(X_j) }{\sum_{\pi\in S_n} \prod_{j=1}^n f_{\theta_{\pi(j)}}(X_j)}. 
\end{align}
This complicated form presents difficulties for both theoretical analysis and computation.

Finally, we recall the following orthogonality relation from \cite{greenshtein2009asymptotic}: 
\begin{align*}
\MSE(\theta,\thetahatS) - \MSE(\theta,\thetahatPI)  &= \bE_\theta\qth{ \| \theta - \thetahatS \|^2 - \| \theta - \thetahatPI \|^2 } \\
&\stepa{=} n\cdot \bE\qth{ (\widetilde{\theta}_1- \bE[\widetilde{\theta}_1|X_1])^2 - (\widetilde{\theta}_1- \bE[\widetilde{\theta}_1|X^n])}^2 \\
&= n\cdot \bE\bqth{ (\bE[\widetilde{\theta}_1|X_1] - \bE[\widetilde{\theta}_1|X^n])^2 \\ &\quad + 2(\bE[\widetilde{\theta}_1|X^n]-\bE[\widetilde{\theta}_1|X_1])\underbrace{(\widetilde{\theta}_1- \bE[\widetilde{\theta}_1|X^n])}_{\bE[\cdot|X^n]=0} } \\
&\stepb{=} n\cdot \bE\bqth{ (\bE[\widetilde{\theta}_1|X_1] - \bE[\widetilde{\theta}_1|X^n])^2}, 
\end{align*}
where (a) uses the postulated Bayes model with the expressions for both oracles, and (b) follows from the tower property of conditional expectation. Consequently,
\begin{align}\label{eq:orthogonality}
\regPI(\widehat{\theta}) - \regS(\widehat{\theta}) &\le \sup_\theta \qth{ \MSE(\theta,\thetahatS) - \MSE(\theta,\thetahatPI) } \nonumber \\
&= n\cdot \sup_\theta \bE\bqth{ (\bE[\widetilde{\theta}_1|X_1] - \bE[\widetilde{\theta}_1|X^n])^2}, 
\end{align}
which serves as our starting point for bounding the difference between EB and compound regrets. 

\subsection{Gaussian model}
In this section we prove \Cref{thm:EB-gaussian} for the Gaussian location model. Using the postulated Bayes model in \Cref{subsec:EB_preliminary} and \eqref{eq:orthogonality}, we only need to show that
\begin{align*}
\bE\qth{\pth{\bE[\widetilde{\theta}_1 | X_1] - \bE[\widetilde{\theta}_1 | X^n]}^2} \le \frac{r_n}{n}. 
\end{align*}
{\blue The proof consists of three parts. First, we apply the noisy interpolation idea of \cite{han2025best} to reduce the difference $\bE[\widetilde{\theta}_1 | X_1] - \bE[\widetilde{\theta}_1 | X^n]$ to $\bE[Z_1 | X_1] - \bE[Z_1 | X^n]$ for a ``noisy'' version $Z$ of $\widetilde{\theta}$. The motivation behind such interpolation is to relate $\bE[Z_1 | X_1] - \bE[Z_1 | X^n]$ to $\KL(P_{Z^n}\|\prod_{i=1}^n P_{Z_i})$ for ``noisy'' permutation mixtures in the second step, via entropic stability arguments. Finally, we invoke our new results \Cref{thm:main_upper} and \Cref{thm:main_upper_dimdep} to upper bound $\KL(P_{Z^n}\|\prod_{i=1}^n P_{Z_i})$.}

We start from the first case, i.e. $|\theta_i|\le h$ for all $i\in [n]$. The first step of the proof is to apply the noisy interpolation idea in \cite{han2025best}: instead of adding $\calN(0,1)$ to $\widetilde{\theta}_i$ for obtaining $X_i$, we first add half the noise $\calN(0,\frac{1}{2})$ to obtain an auxiliary random variable $Z_i$, and then add the remaining half noise $\calN(0,\frac{1}{2})$ to $Z_i$ to arrive at $X_i$. In \cite{han2025best}, it was shown that
\begin{align}\label{eq:Gaussian_identity}
\bE[\widetilde{\theta}_1 | X_1] - \bE[\widetilde{\theta}_1 | X^n] = 2\pth{\bE[Z_1 | X_1] - \bE[Z_1 | X^n]},
\end{align}
so that $\widetilde{\theta}$ in the target inequality can be replaced by $Z$. 
Next we will show that
\begin{align}\label{eq:Tao-Gaussian}
\bE\qth{\pth{\bE[Z_1 | X_1] - \bE[Z_1 | X^n]}^2}  \lesssim \min\left\{
(h^2+1)I(Z_1; X_2^n | X_1),\log n\cdot I(Z_1; X_2^n | X_1) + \frac{1}{n^2}
\right\}. 
\end{align}

The second upper bound of \eqref{eq:Tao-Gaussian} has been shown in \cite{han2025best}; we only need to establish the first bound. To this end, we use the following result on \emph{entropic stability} from \cite{chen2022localization}, which can be viewed as a generalized version of Pinsker's inequality: 
\begin{lemma}\label{lemma:entropic_stability}
Let $\nu$ be a probability measure on $\mathbb{R}$, such that $\var(\calT_t \nu) \le A$ for all $t\in \mathbb{R}$, where $\calT_t \nu(\rmd x) \propto e^{tx}\nu(\rmd x)$ denotes the exponential tilt of $\nu$. Then $\nu$ is $A$-entropically stable in the sense that for every probability measure $\mu$, it holds that
\begin{align*}
\frac{1}{2}\pth{\bE_{\mu}[X] - \bE_{\nu}[X]}^2 \le A\cdot  \KL(\mu\|\nu). 
\end{align*}
\end{lemma}

We apply \Cref{lemma:entropic_stability} to the distribution $\nu = P_{Z_1|X_1}$. Since
\begin{align*}
\nu(\rmd z) &\propto \frac{1}{n}\sum_{i=1}^n \exp\pth{-(z-\theta_i)^2-(X_1-z)^2} \\ &\propto \frac{1}{n}\sum_{i=1}^n \exp\pth{-2\pth{z-\frac{\theta_i+X_1}{4}}^2}
\end{align*}
is a Gaussian mixture with a constant variance for all components, its exponential tilt
\begin{align*}
\calT_t \nu(\rmd z) \propto \frac{1}{n}\sum_{i=1}^n \exp\pth{-2\pth{z-\frac{\theta_i+X_1+t}{4}}^2}
\end{align*}
remains a Gaussian mixture. Since $|\theta_i|\le h$ for all $i\in [n]$, the law of total variance gives $\var(\calT_t \nu)\le \frac{h^2}{16}+\frac{1}{4}$. Therefore, \Cref{lemma:entropic_stability} with $\mu=P_{Z_1|X^n}$ gives
\begin{align*}
\bE\qth{\pth{\bE[Z_1 | X_1] - \bE[Z_1 | X^n]}^2}  &\lesssim (h^2+1) \bE\qth{ \KL(P_{Z_1|X^n} \| P_{Z_1|X_1}) } \\ &= (h^2+1) I(Z_1;X_2^n|X_1),
\end{align*}
which establishes the first inequality of \eqref{eq:Tao-Gaussian}. To proceed, we invoke the arguments in \cite[Lemma 23]{nie2023large} (see also \cite{han2025best}) to obtain
\begin{align}\label{eq:Gaussian-mutualinfo}
    I(Z_1;X_2^n|X_1) \le \frac{1}{n}\KL\pth{ P_{Z^n} \Big\| \prod_{i=1}^n P_{Z_i} } \le \frac{\log(1+\chi^2(\calP,n))}{n}, 
\end{align}
where $\chi^2(\calP,n)$ is the dimension-dependent quantity \eqref{eq:dimension_dep_bound} for the distribution class $\calP = \sth{\calN(\theta_i,\frac{1}{2}): i\in [n]}$. Using the characterization of $\chi^2(\calP,n)$ in \Cref{cor:dim-dep-bound} and distinguishing into three regimes $h\le 1, 1<h\le \sqrt{\log n}$, and $h>\sqrt{\log n}$, from \eqref{eq:Gaussian_identity}--\eqref{eq:Gaussian-mutualinfo} we get
\begin{align*}
\bE\qth{\pth{\bE[\widetilde{\theta}_1 | X_1] - \bE[\widetilde{\theta}_1 | X^n]}^2} \lesssim \frac{h^4 \wedge (h\log^{3/2}n)}{n}. 
\end{align*}
This completes the proof of Case 1 in \Cref{thm:EB-gaussian}.  

Next we turn to Cases 2 and 3 in \Cref{thm:EB-gaussian}. For Case 2, the $s$-subGaussian condition on the empirical distribution $G_n$ ensures that $\|\theta\|_\infty = O_s(\sqrt{\log n})$, so the target bound $r_n = O_s(\log^2 n)$ follows from the $O(h\log^{3/2}n)$ bound in Case 1. For Case 3, the weak $\ell_p$ norm constraint yields
\begin{align*}
\Big|\sth{i\in [n]: |\theta_i| \ge n^{\frac{1}{1+p}}\log^{\frac{1}{2(1+p)}} n} \Big|  \lesssim n\cdot \pth{n^{-\frac{1}{1+p}}\log^{-\frac{1}{2(1+p)}} n}^p = n^{\frac{1}{1+p}}\log^{-\frac{p}{2(1+p)}}n =: K. 
\end{align*}
{\blue Therefore, for $\calP = \sth{\calN(\theta_i,\frac{1}{2}): i\in [n]}$ and $k\ge 2K$, it holds that
\begin{align*}
\Diam \lesssim \frac{(n^{\frac{1}{1+p}}\log^{\frac{1}{2(1+p)}} n)^2}{(k-K)^2} \lesssim \frac{K^2\log n}{k^2}. 
\end{align*}
Consequently, \[\sum_{k=1}^\infty \min\{\Diam, \log n\} \lesssim K\log n + \sum_{k\ge 2K} \frac{K^2\log n}{k^2} \lesssim K\log n.\] 
In addition, \Cref{lemma:chi2-capacity-union} implies that $\Capa\lesssim n^{\frac{1}{1+p}}\log^{\frac{1}{2(1+p)}} n + K \lesssim K\sqrt{\log n}$. By \Cref{thm:main_upper_dimdep}, 
\begin{align*}
\log(1+\chi^2(\calP,n)) &\lesssim \sum_{k=1}^\infty \min\{\Diam, \log n\} + (\Capa+1)\log_+\log \Capa \\
&\lesssim K\log n + K\sqrt{\log n}\log\log n = O\pth{n^{\frac{1}{1+p}} \log^{\frac{2+p}{2(1+p)}} n}, 
\end{align*}
which combined with \eqref{eq:Gaussian_identity}--\eqref{eq:Gaussian-mutualinfo} yields the target upper bound for Case 3. 
}

\subsection{Poisson model}
The proof of \Cref{thm:EB-Poisson} in the Poisson model is conceptually similar but involves additional technical details. We first apply a different interpolation in \cite{han2025best}: condition on $\widetilde{\theta}_i$, we draw $Z_i\sim \Poi(2\widetilde{\theta}_i)$ and $X_i\sim \mathrm{B}(Z_i, \frac{1}{2})$, so that the conditional distribution of $X_i$ given $\widetilde{\theta}_i$ is still $\Poi(\widetilde{\theta}_i)$. The following identity has been shown in \cite{han2025best}: 
\begin{align}\label{eq:Poisson_identity}
\bE[\widetilde{\theta}_1 | X_1] - \bE[\widetilde{\theta}_1 | X^n] = \bE[Z_1 | X_1] - \bE[Z_1 | X^n]. 
\end{align}
Next we show that if $\theta_i\in [0,h]$ for all $i\in [n]$,
\begin{align}\label{eq:Tao-Poisson}
&\bE\qth{\pth{\bE[Z_1 | X_1] - \bE[Z_1 | X^n]}^2}\lesssim \min\left\{
h(h+1) I(Z_1;X_2^n\mid X_1),h\log n\cdot I(Z_1;X_2^n\mid X_1)+\frac{h}{n^2}
\right\}.
\end{align}

To prove the first upper bound of \eqref{eq:Tao-Poisson}, we need the following lemma.

\begin{lemma}\label{lemma:transportation}
Let $\nu = \int \Poi(\theta)\pi(\rmd \theta)$ be a Poisson mixture with $\mathrm{supp}(\pi)\subseteq [0,h]$. Then for any distribution $\mu$ on $\mathbb{R}$ with $0\le \bE_{\mu}[X]\le h$, it holds that
\begin{align*}
\pth{\bE_{\mu}[X] - \bE_{\nu}[X]}^2 \le 2h(h+2)\KL(\mu\|\nu). 
\end{align*}
\end{lemma}
We apply \Cref{lemma:transportation} to the distribution $\nu = P_{Z_1-X_1|X_1}$, with pmf
\begin{align*}
\nu(y) &\propto \frac{1}{n}\sum_{i=1}^n e^{-2\theta_i} \frac{(2\theta_i)^{y+X_1}}{(y+X_1)!} \cdot \binom{y+X_1}{X_1}\frac{1}{2^{y+X_1}}\\ & \propto  \frac{1}{n}\sum_{i=1}^n e^{-2\theta_i} \theta_i^{X_1}\frac{\theta_i^{y}}{y!}, \quad y\in \naturals.
\end{align*}
Therefore, $\nu$ is a Poisson mixture with Poisson rates $\sth{\theta_i: i\in [n]} \subseteq [0,h]$. In addition, 
\begin{align*}
\bE[Z_1-X_1|X^n, \widetilde{\theta}_1] \stepa{=} \bE[Z_1-X_1|X_1, \widetilde{\theta}_1] \stepb{=} \widetilde{\theta}_1 \in [0,h], 
\end{align*}
where (a) is due to the Markov structure $X_2^n - (\widetilde{\theta}_1, X_1) - Z_1$, and (b) uses the fact that conditioned on $\widetilde{\theta}_1$, the random variables $X_1$ and $Z_1- X_1$ are i.i.d. $\Poi(\widetilde{\theta}_1)$. Consequently, the mean of $\mu=P_{Z_1-X_1|X^n}$ is $\bE[Z_1-X_1|X^n]\in [0,h]$ almost surely. Consequently, \Cref{lemma:transportation} implies that
\begin{align*}
    \bE\qth{\pth{\bE[Z_1 | X_1] - \bE[Z_1 | X^n]}^2} &= \bE\qth{\pth{\bE[Z_1 - X_1| X_1] - \bE[Z_1 - X_1 | X^n]}^2} \\
    &\le 2h(h+2) \bE\qth{\KL(P_{Z_1|X^n}\|P_{Z_1|X_1})} \\
    &= 2h(h+2) I(Z_1;X_2^n|X_1), 
\end{align*}
which is the first upper bound of \eqref{eq:Tao-Poisson}. For the second upper bound of \eqref{eq:Tao-Poisson}, note that it is implied by the first upper bound when $h=O(\log n)$; hence in the sequel we assume that $h\ge 64\log n$. Define an event $E=\sth{|Z_1-2X_1|>8\sqrt{h\log n}}$, then Tao's inequality (or Pinsker's inequality) yields
\begin{align}\label{eq:Ec}
&\bE\qth{\pth{\bE[(Z_1 - 2X_1)\indc{E^c}| X_1] - \bE[(Z_1 - 2X_1)\indc{E^c} | X^n]}^2} \nonumber \notag\\
&\le \frac{1}{2}\pth{8\sqrt{h\log n}}^2 \quad \bE\qth{\KL(P_{(Z_1-2X_1)\indc{E^c}|X^n} \| P_{(Z_1-2X_1)\indc{E^c}|X_1} )} \nonumber \notag\\
&= O(h\log n)\cdot I((Z_1-2X_1)\indc{E^c}; X_2^n | X_1) \notag\\
&\le O(h\log n)\cdot I(Z_1; X_2^n| X_1),
\end{align}
where the last step is the data-processing inequality. For the event $E$, note that
\begin{align*}
\bP(E) \le \bP\pth{|Z_1-X_1-\widetilde{\theta}_1| > 4\sqrt{h\log n}} + \bP\pth{|X_1-\widetilde{\theta}_1| > 4\sqrt{h\log n}} = O\pth{\frac{1}{n^4}}, 
\end{align*}
where $Z_1-X_1, X_1\sim \Poi(\widetilde{\theta}_1)$ are conditionally independent given $\widetilde{\theta}_1$, and we have used $\widetilde{\theta}_1\le h$ almost surely and the Poisson tail bound. Consequently, applying Cauchy--Schwarz twice yields
\begin{align}\label{eq:E}
 \bE\qth{\bE[(Z_1 - 2X_1)\indc{E}| X_1]}^2 &\le \bE\qth{(Z_1-2X_1)^2\indc{E}} \nonumber \notag\\
&\le \sqrt{\bE\qth{(Z_1-2X_1)^4} \bP(E)} = O\pth{\frac{h}{n^2}},
\end{align}
where we have used $\bE[(X-Y)^4]=12\theta^2+2\theta$ for independent $X,Y\sim \Poi(\theta)$, and $\widetilde{\theta}_1\le h$ almost surely. The same upper bound also holds for $\bE\qth{\bE[(Z_1 - X_1)\indc{E}| X^n]}^2$. Now using \eqref{eq:Ec}, \eqref{eq:E}, and the triangle inequality arrives at the second upper bound of \eqref{eq:Tao-Poisson}. 

Next, we proceed as in the Gaussian case that \begin{align}\label{eq:Poisson-mutualinfo}
    I(Z_1;X_2^n|X_1) \le \frac{1}{n}\KL\pth{ P_{Z^n} \Big\| \prod_{i=1}^n P_{Z_i} } &\le \frac{\log(1+\chi^2(\calP,n))}{n}, 
\end{align}
where $\chi^2(\calP,n)$ is the dimension-dependent quantity \eqref{eq:dimension_dep_bound} for the distribution class $\calP = \sth{\Poi(2\theta_i): i\in [n]}$. Using the characterization of $\chi^2(\calP,n)$ in \Cref{cor:dim-dep-bound} and distinguishing into three regimes $h\le 1, 1<h\le \log n$, and $h>\log n$, from \eqref{eq:Poisson_identity}--\eqref{eq:Poisson-mutualinfo} we get
\begin{align*}
\bE\qth{\pth{\bE[\widetilde{\theta}_1 | X_1] - \bE[\widetilde{\theta}_1 | X^n]}^2} \lesssim \frac{h^3 \wedge (h\log n)^{3/2}}{n}, 
\end{align*}
which is the desired bound for Case 1. For Case 2, the subexponential condition on the empirical distribution $G_n$ implies that $\|\theta\|_\infty = O_s(\log n)$, so using the bound for Case 1 with $h=O_s(\log n)$ yields the target bound $O(\log^3 n)$. 

Finally we handle Case 3 with a bounded $p$-th moment. We first modify the definition of $E$ by choosing
\begin{align}\label{eq:h}
h \asymp n^{\frac{2}{2p+1}}\log^{\frac{1}{2p+1}} n. 
\end{align}
For this specific choice of $h$, the inequality \eqref{eq:Ec} still holds, but \eqref{eq:E} now becomes
\begin{align}\label{eq:E_smallh}
\bE\qth{\bE[(Z_1 - 2X_1)\indc{E, \widetilde{\theta}_1\le h}| X_1]}^2 = O\pth{\frac{h}{n^2}}. 
\end{align}
For the remaining case $\widetilde{\theta}_1>h$, we simply bound
\begin{align}\label{eq:E_largeh}
\bE\qth{\bE[(Z_1 - 2X_1)\indc{E, \widetilde{\theta}_1> h}| X_1]}^2 &\le \bE[(Z_1-2X_1)^2 \indc{E, \widetilde{\theta}_1> h}] \nonumber \\
&\le \bE[(Z_1-2X_1)^2 \indc{\widetilde{\theta}_1> h}] \nonumber \\
&= \frac{1}{n}\sum_{i\in [n]: \theta_i > h} \bE[(Z_1-2X_1)^2 | \widetilde{\theta}_1 = \theta_i] \nonumber \\
&= \frac{2}{n}\sum_{i\in [n]: \theta_i > h} \theta_i \nonumber \\
&\le \frac{2}{n}\sum_{i\in [n]: \theta_i > h} \frac{\theta_i^p}{h^{p-1}} = O\pth{\frac{1}{h^{p-1}}}, 
\end{align}
where the last step uses the moment bound of $G_n$. Therefore, by \eqref{eq:Poisson_identity}, \eqref{eq:Ec}, \eqref{eq:Poisson-mutualinfo}, \eqref{eq:E_smallh}, and \eqref{eq:E_largeh}, we arrive at
\begin{align}\label{eq:regret_intermediate}
&n\bE\qth{\pth{\bE[\widetilde{\theta}_1 | X_1] - \bE[\widetilde{\theta}_1 | X^n]}^2} \lesssim (h\log n)\log(1+\chi^2(\calP,n)) + \frac{h}{n} + \frac{n}{h^{p-1}}. 
\end{align}
Finally we use \Cref{thm:main_upper_dimdep} to upper bound the dimension-dependent quantity $\chi^2(\calP,n)$. {\blue By Markov's inequality applied to the bounded $p$-th moment, 
\begin{align*}
|\sth{i\in [n]: \theta_i\ge h}| \lesssim \frac{n}{h^p} =: K. 
\end{align*}
Therefore, for $\calP=\sth{\Poi(2\theta_i):i\in [n]}$ and $k\ge 2K$, 
\begin{align*}
\Diam \lesssim \frac{h}{(k-K)^2} \lesssim \frac{h}{k^2}, 
\end{align*}
so that 
\begin{align*}
\sum_{k=1}^\infty \min\sth{\Diam,\log n}  &\lesssim K\log n + \sum_{k\ge 2K} \frac{h}{k^2} \\ & \lesssim \frac{n\log n}{h^p} + \frac{h^{2p+1}}{n} \\& \overset{\prettyref{eq:h}}{=} \frac{n\log n}{h^p} + \sqrt{h\log n}. 
\end{align*}
In addition, \Cref{lemma:chi2-capacity-union} implies that $\Capa\lesssim \sqrt{h} + K = \sqrt{h}+\frac{n}{h^p}$. By \Cref{thm:main_upper_dimdep}, 
\begin{align}\label{eq:chi-squared-bound}
\log(1+\chi^2(\calP,n)) &\lesssim \sum_{k=1}^\infty \min\{\Diam, \log n\} + (\Capa+1)\log_+\log \Capa \notag\\ &= O\pth{\sqrt{h\log n} + \frac{n\log n}{h^p}}. 
\end{align}
} The desired upper bound then follows from \eqref{eq:regret_intermediate} and \eqref{eq:chi-squared-bound}, with the choice of $h$ given in \eqref{eq:h}.

%% file: Discussion.tex
\section{Discussion}\label{sec:discussions}
In this paper, we have advanced the understanding of the information geometry underlying permutation mixtures. Our first contribution is an improved mean-field analysis in both infinite-dimensional and finite-dimensional settings, achieved by establishing a link between the eigenvalues of the channel overlap matrix and the R\'enyi partition diameter, a geometric property of the channel $\mathcal{P}$. These refinements close several polynomial gaps left open in \cite{han2024approximate}, and are precisely sharp enough to demonstrate the tightness of approximating PI decision rules by separable ones, thereby establishing strong equivalence between two notions of compound regret in several non-asymptotic regimes.

We end this paper with several discussions and open questions. 

\paragraph{Removing extra log factors} Our upper bound in \Cref{thm:main_upper} involves a $\log\log$ factor on the $\chi^2$ channel capacity, which was not present in \cite{han2024approximate}. Technically, this is due to the logarithmic factor in the higher-order Cheeger inequality. A second logarithmic gap occurs when $\Diam \asymp k^{-1}$ (e.g. in the normal variance model of \Cref{cor:dim-indep-bound}), where the harmonic sum $\sum_k \Diam$ contributes a logarithmic factor to the upper bound, whereas the lower bound essentially depends only on $\max_k k \Diam$. We conjecture that the upper bound can be improved, as the current link by Cheeger inequality upper bounds $\prod_{k\ge 2} \max_A  \frac{1}{1-\lambda_k(A)}$ rather than $ \max_A  \prod_{k\ge 2} \frac{1}{1-\lambda_k(A)}$. 

\paragraph{Bridging different upper bounds}Our upper bound in \Cref{thm:main_upper} for $\Capa \ge 1$ appears disconnected from the $O(\Capa^2)$ bound known for $\Capa < 1$. In fact, {\blue a recent work \cite{han2026elementary} shows the following strengthening of \eqref{eq:A-upper}:}
\begin{align*}
\chi^2(\bP_n \| \bQ_n) \le C\bpth{ \prod_{i=2}^n \frac{1}{1-\lambda_i(A)} - 1 - \sum_{i=2}^n \lambda_i(A)}.  
\end{align*}
This bound recovers both \eqref{eq:A-upper} and the $O(\Capa^2)$ upper bound in the small eigenvalue limit. Moreover, it narrows the gap with the lower bound in \Cref{thm:main_lower}, whose first-order term in the Taylor expansion is quadratic in eigenvalues. By linking eigenvalues to R\'enyi partition diameters, the resulting bound unifies the two regimes covered by existing upper bounds.

\paragraph{Improving regret lower bounds} The possible logarithmic gap between $\regS$ and $\regPI$ in Case~3 of \Cref{cor:EB-gaussian} and \Cref{cor:EB-Poisson} may be partly attributable to a potentially loose lower bound on $\regS$. In particular, both \cite{shen2022empirical} and \Cref{append:regret_LB} derive lower bounds by essentially counting the number of ``local intervals'', whereas \cite{kim2014minimax,polyanskiy2021sharp} show that the regret is already at least logarithmic even for a single local interval. This suggests that combining these techniques could yield a sharper lower bound on $\regS$, potentially leading to $\regPI \asymp \regS$.

\paragraph{Generalization to other EB models} The current mean-field analysis of the PI oracle critically relies on a divisible structure of the underlying distribution, which is in sharp contrast to the mean-field analysis of permutation mixtures where no such structure is required. Nevertheless, this approach successfully unifies the upper bounds of \cite{greenshtein2009asymptotic} and \cite{han2025best}, whose underlying ideas appear quite different. An interesting direction is to better understand the scope and limitations of this analysis, and to generalize it to EB models (such as the normal variance model) that lack divisibility structures.

\paragraph{Equivalence between EB and compound regrets} Beyond the equivalence between the two notions of compound regret, it is natural to ask whether there exists a black-box reduction from the EB regret to the compound regret for a given estimator. For instance, \cite{polyanskiy2021sharp} noted that it remains open whether Robbins's estimator can also achieve a logarithmic compound regret. The framework of permutation mixtures appears relevant here, as the marginal distributions of $X^n$ coincide with $\bQ_n$ and $\bP_n$ in the EB and compound settings, respectively. The main challenge lies in the dependence of the estimation risk on $\theta^n$, whose distributions differ substantially between the two settings.

\section*{Acknowledgement} We thank Yury Polyanskiy for posing the question of the equivalence between two notions of compound regrets in Gaussian location models with a subGaussian prior, which served as the primary motivation for this work. We thank Jonathan Niles-Weed for drawing our attention to the technique of entropic stability used in \Cref{lemma:entropic_stability}, and Cun-Hui Zhang for helpful discussions on the separable and PI oracles. {\blue We also thank anonymous reviewers for their helpful comments, and in particular, for suggesting the connection between \Cref{thm:main_upper} and an integral of covering numbers.}

%% file: Deferred_Proofs.tex
\section{Deferred Proofs}
\subsection{Proof of \Cref{cor:dim-indep-bound}}
We begin with two useful lemmas to upper bound the $\chi^2$ channel capacity $\Capa$. 

\begin{lemma}[Lemma C.1 of \cite{han2024approximate}]\label{lemma:chi2-capacity-union}
Let $\calP_1, \ldots, \calP_m$ be families of probability distributions over the same space. Then
\begin{align*}
\sfC_{\chi^2}\pth{\bigcup_{i=1}^m \calP_i} \le \sum_{i=1}^m \sfC_{\chi^2}\pth{\calP_i} + m - 1.
\end{align*}
\end{lemma}

\begin{lemma}\label{lemma:chi2-capacity-ratio}
Suppose there is a probability distribution $Q$ (not necessarily in $\calP$) such that
\begin{align*}
\sup_{P\in \calP} \Big\| \frac{\rmd P}{\rmd Q} \Big\|_\infty \le A. 
\end{align*}
Then $\Capa\le A-1$. 
\end{lemma}
\begin{proof}
Let $(P_\theta)$ be a parametrization of $\calP$, and $\rho$ be any prior distribution of $\theta$. Then
\begin{align*}
\bE_{\theta\sim \rho}\qth{\chi^2(P_\theta\| \bE_{\theta'\sim\rho}[P_{\theta'}] )} &= \bE_{\theta\sim \rho}\qth{ \int \frac{(\rmd P_{\theta})^2}{\bE_{\theta'\sim \rho}[\rmd P_{\theta'}]} } - 1 \\
&= \int \frac{\bE_{\theta\sim\rho}[(\rmd P_{\theta})^2] }{\bE_{\theta\sim\rho}[(\rmd P_{\theta})]} - 1 \\
&\le \int \frac{\bE_{\theta\sim\rho}[(\rmd P_{\theta})(A\rmd Q)] }{\bE_{\theta\sim\rho}[(\rmd P_{\theta})]} - 1  \\ &= A-1, 
\end{align*}
where the inequality follows from the deterministic relation $\rmd P_\theta\le A\rmd Q$. 
\end{proof}

Also, we recall from \cite[Lemma 6.1]{han2024approximate} (or equivalently, \Cref{thm:main_lower} with $n=2$) a weaker lower bound:
\begin{align}\label{eq:LB_existing}
\chi^2(\calP) + 1 = \exp(\Omega(\sfD_1(\calP))). 
\end{align}
Now we are ready to prove \Cref{cor:dim-indep-bound}. 

\paragraph{Normal mean family} Since
\begin{align*}
D_{\frac{1}{2}}(\calN(\theta_1,1), \calN(\theta_2,1)) = \frac{(\theta_1-\theta_2)^2}{4}, 
\end{align*}
by partitioning $[-\mu,\mu]$ into $k$ subintervals of equal length we obtain $\sfD_k(\calP) \lesssim \mu^2/k^2$. In addition, \cite[Section C.2]{han2024approximate} shows that $\Capa \lesssim \mu$. Since $\sum_{k\ge 1} \frac{1}{k^2}<\infty$, \Cref{thm:main_upper} gives the upper bound $\chi^2(\calP) = \exp(O(\mu^2 + \mu\log\log \mu)) = \exp(O(\mu^2))$. The claimed lower bound directly follows from \eqref{eq:LB_existing} and $\sfD_1(\calP)=\mu^2$. 

\paragraph{Poisson family} Let $\Delta:=\sqrt{M}-\sqrt{m}\ge 1$. Since
\begin{align*}
D_{\frac{1}{2}}(\Poi(\lambda_1), \Poi(\lambda_2)) = (\sqrt{\lambda_1}-\sqrt{\lambda_2})^2, 
\end{align*}
by considering the uniform partition on $\sqrt{\lambda}\in [\sqrt{m}, \sqrt{M}]$, we obtain $\sfD_k(\calP)\lesssim \Delta^2/k^2$. In addition, \cite[Section C.2]{han2024approximate} shows that $\Capa \lesssim \Delta$. Therefore, \Cref{thm:main_upper} gives the upper bound $\chi^2(\calP) = \exp(O(\Delta^2 + \Delta\log\log \Delta)) = \exp(O(\Delta^2))$. The lower bound follows from \eqref{eq:LB_existing} and $\sfD_1(\calP)=\Delta^2$. 

\paragraph{Normal variance family} Since
\begin{align*}
D_{\frac{1}{2}}(\calN(0,\sigma_1^2), \calN(0,\sigma_2^2)) = \log\pth{\frac{\sigma_1^2+\sigma_2^2}{2\sigma_1\sigma_2}}, 
\end{align*}
by choosing $\sigma_1=e^{-M}$ and $\sigma_2 = e^M$, it is clear that $\sfD_1(\calP) \asymp M$. In addition, by considering a geometric partition $[e^{-M},e^{-M+2M/k}), [e^{-M+2M/k},e^{-M+4M/k}), \dots$ $, [e^{M-2M/k},e^{M})$ for the domain of $\sigma\in [e^{-M}, e^M]$, we obtain $\sfD_k(\calP)\lesssim M/k$ for $k=O(M)$. As for $\Capa$, we again consider a geometric partition $[e^{-M},\rho e^{-M}), [\rho e^{-M}, \rho^2 e^{-M}), \dots, [\rho^{s-1} e^{-M}, \rho^{s} e^{-M}]$ with a fixed constant $\rho\in (1,\sqrt{2})$, and $\rho^{s}e^{-M}=e^M$. Since $1<\rho<\sqrt{2}$ and
\begin{align*}
\chi^2(\calN(0,\sigma_1^2) \| \calN(0,\sigma_2^2)) = \frac{\sigma_2^2}{\sqrt{\sigma_1^2(2\sigma_2^2-\sigma_1^2)}}-1, \\ \quad \text{if } \sigma_1^2 < 2\sigma_2^2,
\end{align*}
the $\chi^2$ diameter for each subset in the partition is $O(1)$; \Cref{lemma:chi2-capacity-union} then gives $\Capa=O(s)=O(M)$. Consequently, \Cref{thm:main_upper} gives an upper bound
\begin{align*}
\log(1+\chi^2(\calP)) &= O\bpth{ \sum_{k=1}^{O(M)}\frac{M}{k} + M\log\log M } \\ &= O(M\log M). 
\end{align*}
The lower bound again follows from \eqref{eq:LB_existing} and $\sfD_1(\calP)\asymp M$.

\paragraph{Normal mean family in dimension $d\ge 2$} Again, 
\begin{align*}
D_{\frac{1}{2}}(\calN(\theta_1,I_d), \calN(\theta_2,I_d)) = \frac{\|\theta_1-\theta_2\|^2}{4}. 
\end{align*}
By standard covering arguments, the $d$-dimensional ball with radius $r$ can be covered by $O_d(r^d)$ unit balls. By applying such a cover with $k=O_d(r^d)$ and proper scaling, it is clear that $\sfD_k(\calP)\lesssim \mu^2/k^{2/d}$. In addition, in view of
\begin{align*}
\chi^2(\calN(\theta_1,I_d)\|\calN(\theta_2,I_d)) = e^{\|\theta_1-\theta_2\|^2}-1, 
\end{align*}
the $\chi^2$ diameter of each unit ball is $O(1)$. Therefore, \Cref{lemma:chi2-capacity-union} implies that $\Capa=O_d(\mu^d)$. Finally, by \Cref{thm:main_upper}, we have the upper bound
\begin{align*}
\log(1+\chi^2(\calP)) &= O_d\pth{\sum_{k=1}^{O_d(\mu^d)} \frac{\mu^2}{k^{2/d}} + \mu^d\log\log \mu} \\ &= \begin{cases}
    O(\mu^2\log \mu) & \text{if } d = 2, \\
    O_d(\mu^d\log\log\mu) & \text{if } d \ge 3.
\end{cases}
\end{align*}

Next we prove the lower bound. WLOG assume that $\mu\ge C_0(d)$ is large enough, where $C_0(d)>0$ is a large constant depending on $d$ and to be chosen later; in fact, if $\mu<C_0(d)$, we simply invoke the lower bound \eqref{eq:LB_existing} with $\sfD_1(\calP)= \mu^2$. Next consider the packing
\begin{align*}
\sth{\theta_1,\dots,\theta_M} &:= \sth{-\frac{\mu}{\sqrt{d}}, -\frac{\mu}{\sqrt{d}} + C_1, -\frac{\mu}{\sqrt{d}} + 2C_1, \dots, \frac{\mu}{\sqrt{d}} }^d \\ &\subseteq \sth{\theta\in \bR^d: \|\theta\|\le \mu }, 
\end{align*}
where $C_1=C_1(d)>0$ is another large constant depending only on $d$. Note that $M=\Omega_d(\mu^d)$, and we can choose $C_0(d)$ to be a large enough constant depending on $C_1$ and $d$ such that $2\mu/\sqrt{d}$ is an integral multiple of $C_1$. Now let $P_i = \calN(\theta_i,I_d)\in \calP$ for $i\in [M]$, and construct the channel overlap matrix $A$ (see \Cref{defn:channel-overlap-matrix}) based on $P_1,\dots,P_M$. Note that for any $i\in [M]$,
\begin{align*}
\sum_{j\neq i}A_{ij} &= \sum_{j\neq i} \int \frac{\rmd P_i \rmd P_j}{\sum_{k=1}^M \rmd P_k} \\& \le \sum_{j\neq i} \int \frac{\rmd P_i \rmd P_j}{\rmd P_i + \rmd P_j} \\ &\le \frac{1}{2}\sum_{j\neq i}\int \sqrt{\rmd P_i \rmd P_j} \\&= \frac{1}{2}\sum_{j\neq i} \exp\pth{-\frac{\|\theta_i-\theta_j\|^2}{8}}.
\end{align*}
In addition, by the grid structure of $\sth{\theta_1,\dots,\theta_M}$, we have
\begin{align*}
&\sum_{j\neq i} \exp\pth{-\frac{\|\theta_i-\theta_j\|^2}{8}} \\&= \sum_{j=1}^M \exp\pth{-\frac{\|\theta_i-\theta_j\|^2}{8}} - 1 \\
&\le \qth{1+2\exp\pth{-\frac{C_1^2}{8}}+2\exp\pth{-\frac{(2C_1)^2}{8}}+\cdots}^d - 1,
\end{align*}
which is smaller than $1$ for $C_1=C_1(d)$ large enough. Combining the above two inequalities, for all $i\in [M]$ we have
\begin{align*}
A_{ii} = 1-\sum_{j\neq i} A_{ij} \ge \frac{1}{2},
\end{align*}
so that the second inequality of \Cref{thm:main_lower} gives that $\chi^2(\calP)\ge \frac{1}{\sqrt{M}}\exp(\Omega(M))-1=\exp(\Omega_d(\mu^d))$ for all $d\ge 2$, as claimed. 

\paragraph{Discrete distribution family} Since for $P_1,P_2\in \calP$, 
\begin{align*}
\int \sqrt{\rmd P_1 \rmd P_2} = \sum_{i=1}^m \sqrt{P_1(i) P_2(i)}\ge \sqrt{P_1(1)P_2(1)} \ge \varepsilon, 
\end{align*}
we conclude that $\sfD_1(\calP)\le 2\log(1/\varepsilon)$. In addition, by choosing $Q=\Unif([m])$ in \Cref{lemma:chi2-capacity-ratio}, we obtain $\Capa\le m-1$. Consequently, the upper bound in \cite{han2024approximate} gives
\begin{align*}
\log(\chi^2(\calP)+1) = O\pth{ (\sfD_1(\calP)+1)\Capa } = O\pth{m\log\frac{1}{\varepsilon}}. 
\end{align*}

For the lower bound, if $m\ge 3$, we apply \Cref{thm:main_lower} to the following distributions $P_2,\dots,P_m\in \calP$: 
\begin{align*}
P_i = \varepsilon \delta_1 + (1-\varepsilon) \delta_i, \quad 2\le i\le m.
\end{align*}
By simple algebra, the channel overlap matrix $A\in \mathbb{R}^{(m-1)\times (m-1)}$ is
\begin{align*}
A = \frac{\varepsilon}{m-1} J_{m-1} + (1-\varepsilon) I_{m-1}, 
\end{align*}
where $I, J$ are the identity and all-ones matrices, respectively. It is then clear that the eigenvalues of $A$ are $1$ (with multiplicity $1$) and $1-\varepsilon$ (with multiplicity $m-2$), so that \Cref{thm:main_lower} gives
\begin{align*}
\chi^2(\calP) \ge \pth{\frac{1}{1-(1-\varepsilon)^2}}^{\frac{m-2}{2}} = \exp\pth{\Omega\pth{m\log\frac{1}{\varepsilon}}}, 
\end{align*}
thanks to the assumptions $\varepsilon<1/2$ and $m\ge 3$. For $m=2$, we pick $P_1=(\varepsilon,1-\varepsilon)$ and $P_2 = (1-\varepsilon,\varepsilon)$, so that the channel overlap matrix is $A=2\varepsilon(1-\varepsilon)J_2 + (1-2\varepsilon)^2 I_2$. Therefore, $\lambda_2(A) = (1-2\varepsilon)^2$, and \Cref{thm:main_lower} again gives
\begin{align*}
\chi^2(\calP) \ge \frac{1}{\sqrt{1-(1-2\varepsilon)^4}} = \exp\pth{\Omega\pth{\log \frac{1}{\varepsilon}}}. 
\end{align*}

\subsection{Proof of \Cref{cor:dim-dep-bound}}
\paragraph{Normal mean family} {\blue The first upper bound $O(\mu^4)$ for $\mu\le 1$ is a dimension-independent bound established in \cite[Corollary 1.3]{han2024approximate}. For $\mu \ge 1$, since $\Capa\asymp \mu$ and $\Diam\asymp \frac{\mu^2}{k^2}$, the rest of the upper bounds follow from 
\begin{align*}
   \sum_{k=1}^\infty \min\{\Diam, \log n\} \asymp \min\{\mu^2, \mu\sqrt{\log n}\}, 
\end{align*}
which dominates $\mu\log_+\log \mu$ for $1\le \mu\le n^{0.99}$. 
}

Next we establish the lower bound by constructing specific instances $P_1,\dots,P_n$. When $\mu\le \sqrt{\log n}$, we choose $P_1=\dots=P_{m}=\calN(-\mu,1)$, and $P_{m+1}=\dots=P_n=\calN(\mu,1)$, where WLOG we assume that $n=2m$ is an even integer since $n\ge 2$ and $\chi^2(\calP,n)$ is non-decreasing in $n$. This is the same toy example in \cite{han2024approximate}, for which \cite[Sec. 4.1]{han2024approximate} shows the following expression for $\chi^2(\bP_n \| \bQ_n)$:
\begin{align}\label{eq:chi-squared-expression}
\chi^2(\bP_n \| \bQ_n) = \sum_{\ell=1}^{m} \frac{\binom{m}{\ell}^2}{\binom{2m}{2\ell}} f(\mu)^{2\ell}, 
\end{align}
where
\begin{align*}
f(\mu) &= \exp\pth{-\frac{\mu^2}{2}}\bE_{Z\sim \calN(0,1)}\qth{\frac{\sinh^2(\mu Z)}{\cosh(\mu Z)}} \\
&= 1 - \exp\pth{-\frac{\mu^2}{2}} \bE_{Z\sim \calN(0,1)}\qth{\frac{1}{\cosh(\mu Z)}} \\&\ge 1-\exp\pth{-\frac{\mu^2}{2}}. 
\end{align*}
By Stirling's approximation $\sqrt{2\pi n}(\frac{n}{e})^n\le n!\le 2\sqrt{2\pi n}(\frac{n}{e})^n$, we have the lower bound
\begin{align*}
\frac{\binom{m}{\ell}^2}{\binom{2m}{2\ell}} \ge \frac{1}{8\sqrt{\pi}} \sqrt{\frac{m}{\ell(m-\ell)}} \ge \frac{1}{8\sqrt{\pi \ell}}. 
\end{align*}
Plugging these lower bounds into \eqref{eq:chi-squared-expression} yields
\begin{align*}
\chi^2(\bP_n \| \bQ_n) \ge  \sum_{\ell=1}^m \frac{1}{8\sqrt{\pi \ell }}\pth{1-e^{-\mu^2/2}}^{2\ell}. 
\end{align*}
When $\mu\le 1$, the contribution of $\ell=1$ gives the target lower bound $\Omega(\mu^4)$. When $1<\mu\le \sqrt{\log n}$, we have
\begin{align*}
\chi^2(\bP_n \| \bQ_n) &\ge   \sum_{\ell=1}^m \frac{1}{8\sqrt{\pi\ell}}\pth{1-e^{-\mu^2/2}}^{2\ell} \\&\ge \sum_{\ell=1}^{\lfloor e^{\mu^2/2} \rfloor}\frac{1}{8\sqrt{\pi\ell}} \pth{1-e^{-\mu^2/2}}^{2\ell} \\
&\ge \frac{1}{8\sqrt{\pi \exp(\mu^2/2)}} \sum_{\ell=1}^{\lfloor e^{\mu^2/2} \rfloor} \pth{1-e^{-\mu^2/2}}^{2\ell} \\
&= \frac{(1-e^{-\mu^2/2})^2}{8\sqrt{\pi \exp(\mu^2/2)}} \cdot \frac{1-(1-e^{-\mu^2/2})^{2\lfloor e^{\mu^2/2} \rfloor}}{1-(1-e^{-\mu^2/2})^2} \\&= \Omega\pth{e^{\mu^2/4}}. 
\end{align*}
These prove the first two lower bounds. 

The lower bound for $\mu>\sqrt{\log n}$ requires a different construction. First, since the case $\mu\asymp \sqrt{\log n}$ can be covered by the previous $\exp(\Omega(\log n))$ lower bound (by making $\mu$ smaller if necessary), in the sequel we assume that $\mu>\sqrt{C_0\log n}$ for a large absolute constant $C_0>0$. Next we partition $[-\mu,\mu]$ into $m=\frac{\mu}{2\sqrt{\log n}}$ intervals $I_1,\dots,I_m$ of the same length $4\sqrt{\log n}$. Since $\sqrt{C_0\log n}<\mu\le n^{0.99}$ with a large constant $C_0$, WLOG we assume that both $m$ and $\frac{n}{m}$ are integers. Let $\mu_1,\dots,\mu_m$ be the midpoints of $I_1,\dots,I_m$, and construct
\begin{align*}
    P_i = \calN(\mu_{\lceil im/n \rceil}, 1), \qquad i\in [n].  
\end{align*}
In other words, each of $\mu_1,\dots,\mu_m$ appears $\frac{n}{m}$ times in the normal means. 

To lower bound $\chi^2(\bP_n \| \bQ_n)$, we use the data processing inequality to write
\begin{align*}
\chi^2(\bP_n \| \bQ_n) \ge \sup_E \frac{(\bP_n(E) - \bQ_n(E))^2}{\bQ_n(E)}, 
\end{align*}
where the supremum is taken over all possible events $E$. To the desired lower bound, it remains to construct an event $E$ such that
\begin{align}\label{eq:E_Gaussian}
\bP_n(E) = \Omega(1), \qquad \bQ_n(E) = \exp\pth{-\Omega\pth{\mu\sqrt{\log n}}}. 
\end{align}
To show \eqref{eq:E_Gaussian}, we define
\begin{align*}
E = \sth{ |\{i\in [n]: X_i\in I_j\}| = \frac{n}{m} \text{ for all }j\in [m]},
\end{align*}
where we recall that $X_i\sim \calN(\theta_i,1)$, with $(\theta_1,\dots,\theta_n)$ drawn either without replacement or with replacement from the normal means. In other words, the event $E$ states that similar to the normal means $\theta_1,\dots,\theta_n$, the observations $X_1,\dots,X_n$ are also evenly distributed in the intervals $I_1,\dots,I_m$. 

To lower bound $\bP_n(E)$, note that $E$ holds if $|X_i - \theta_{i}|<2\sqrt{\log n}$ for all $i\in [n]$. By Gaussian tail bound, 
\begin{align*}
\bP\pth{ |X_i - \theta_{i}| \ge 2\sqrt{\log n} } \le 2\exp\pth{-\frac{1}{2}\bpth{2\sqrt{\log n}}^2} = \frac{2}{n^2}. 
\end{align*}
Therefore, for $C_1>8$, the union bound over $i\in [n]$ yields $\bP_n(E) \ge 1-\frac{2}{n}= \Omega(1)$. 

To upper bound $\bQ_n(E)$, we need the help from an additional event: 
\begin{align*}
E' = \sth{ \sum_{i=1}^n \indc{|X_i - \theta_{i}|\ge 2\sqrt{\log n}} \le m}.
\end{align*}
In other words, the event $E'$ states that the number of observations belonging to a different interval from their means is at most $m$. By binomial tail bound, it is clear that
\begin{align*}
\bQ_n(E'^c) \le \bP\pth{ \mathrm{B}(n,2n^{-C_1/8})> m} = n^{-\Omega((\frac{C_1}{8}-1)\cdot m)},
\end{align*}
where with a slight abuse of notation, we extend the definition of $\bQ_n$ to the joint distribution of $(\theta^n, X^n)$. Since the event $E\cap E'$ implies that $\theta_1,\dots,\theta_n$ become evenly distributed in the intervals $I_1,\dots,I_m$ after changing at most $m$ of them, we have
\begin{align*}
\bQ_n(E) &\le \bQ_n(E'^c) + \bQ_n(E\cap E') \\
&\le n^{-\Omega((\frac{C_1}{8}-1)\cdot cm)} \\ &\quad+ \bP\pth{ \| \mathrm{Multi}(n; \Unif([m])) - (\frac{n}{m}, \dots, \frac{n}{m}) \|_1 \le m}. 
\end{align*}
Here $\mathrm{Multi}(n; p)$ denotes the multinomial distribution with $n$ draws from the pmf $p=(p_1,\dots,p_m)$. To upper bound the final multinomial probability, note that the number of integer vectors $(x_1,\dots,x_m)$ with $\ell_1$ norm at most $t$ is
\begin{align*}
\le 2^m \cdot \Big|\bsth{y_1,\dots,y_m\in \naturals: \sum_{j=1}^m y_j \le t} \Big| &= 2^m \binom{t+m}{m} \le 2^{2m+t}.
\end{align*}
Moreover, the multinomial pmf is at most 
\begin{align*}
\frac{1}{m^n}\binom{n}{\frac{n}{m}, \dots, \frac{n}{m}} \le \frac{2\sqrt{2\pi n}}{(2\pi \frac{n}{m})^{\frac{m}{2}}} \le \exp(-c_0m\log n)
\end{align*}
by Stirling's approximation and the assumption $m\le \mu\le n^{0.99}$, with an absolute constant $c_0>0$. Consequently, 
\begin{align*}
&\bP\pth{ \| \mathrm{Multi}(n; \Unif([m])) - (\frac{n}{m}, \dots, \frac{n}{m}) \|_1 \le m} \\ &\le 2^{3m} \exp(-c_0m\log n) = \exp(-\Omega(m\log n)).
\end{align*}

A combination of the above inequalities leads to the desired bound
\begin{align*}
    \bQ_n(E) = \exp\pth{-\Omega(m\log n)}=\exp\pth{-\Omega\bpth{\mu\sqrt{\log n}}}, 
\end{align*}
which establishes \eqref{eq:E_Gaussian}.

\paragraph{Poisson family} The analysis in the Poisson family is entirely similar. For the upper bounds, the $O(M^2)$ upper bound for $M\le 1$ follows from \cite[Corollary 1.3]{han2024approximate}, and {\blue 
the remaining upper bounds follow from
\begin{align*}
\sum_{k=1}^\infty \min\{\Diam, \log n\} &\asymp \sum_{k=1}^\infty \min\sth{\frac{M}{k^2}, \log n} \\& \asymp \min\{M, \sqrt{M\log n}\}, 
\end{align*}
which dominates $\sqrt{M}\log_+\log M$ for $1\le M\le n^{1.99}$. 
} 

For the lower bound with $M\le \log n$, we again construct a two-component example $P_1=\dots=P_m=\Poi(0)$ and $P_{m+1}=\dots=P_n=\Poi(M)$, with $n=2m$. An argument entirely similar to \eqref{eq:chi-squared-expression} yields
\begin{align*}
\chi^2(\bP_n \| \bQ_n) = \sum_{\ell=1}^{m} \frac{\binom{m}{\ell}^2}{\binom{2m}{2\ell}} f(M)^{2\ell}, 
\end{align*}
with
\begin{align*}
f(M) &= \frac{1}{2}\chi^2\bpth{\Poi(0) \Big\| \frac{\Poi(0)+\Poi(M)}{2}} \\ & \quad+ \frac{1}{2}\chi^2\bpth{\Poi(M) \Big\| \frac{\Poi(0)+\Poi(M)}{2}} \\
&= \frac{1}{2}\pth{ \frac{2(1+e^{-2M})}{1+e^{-M}} + \sum_{k=1}^\infty 2\bP(\Poi(M)=k) } - 1 \\ &= \frac{1-e^{-M}}{1+e^{-M}}. 
\end{align*}
Now a similar analysis to the Gaussian case yields $\chi^2(\bP_n\|\bQ_n) = \Omega(M^2)$ if $M\le 1$, and $\chi^2(\bP_n\|\bQ_n) = e^{\Omega(M)}$ if $1<M\le \log n$. For $M>\log n$, we modify the Gaussian construction motivated by the Poisson tail bound: the intervals $I_1,\dots,I_m$ are now chosen to be
\begin{align*}
I_j = \qth{ 10(j-1)^2\log n, 10j^2\log n }, \quad \text{ with } m \asymp \sqrt{\frac{M}{\log n}}. 
\end{align*}
These intervals are constructed so that if $x_j$ is the midpoint of $I_j$, then $\bP(\Poi(x_j)\notin I_j)\le \frac{2}{n^2}$ by the Poisson tail bound. Therefore, an entirely similar analysis in the Gaussian case can be carried out to conclude that
\begin{align*}
\chi^2(\bP_n \| \bQ_n) = \exp\pth{\Omega(m\log n)}=\exp\pth{\Omega\bpth{\sqrt{M\log n}}}, 
\end{align*}
where the required condition $m\le n^{0.999}$ is ensured by the assumption $M\le n^{1.99}$. 

\subsection{Proof of \Cref{lemma:transportation}}
The moment generating function (MGF) of $\nu$ is
\begin{align*}
\bE_{\nu} [e^{\lambda X}] &= \bE_{\theta\sim \pi}\qth{\exp\pth{\theta(e^{\lambda}-1)}} \\ &\le \exp\pth{\bE_{\theta\sim \pi}[\theta](e^{\lambda}-1) + \frac{(e^\lambda-1)^2 h^2}{8}},
\end{align*}
where the last inequality is Hoeffding's upper bound on the MGF of bounded random variables. By Donsker--Varadhan \cite[Theorem 4.6]{polyanskiy2024information}, for every $\lambda\in \mathbb{R}$ we have
\begin{align*}
\KL(\mu\|\nu) &\ge \lambda \bE_{\mu}\qth{X} - \log \bE_{\nu} [e^{\lambda X}] \\
&\ge \lambda \bE_{\mu}\qth{X} - \pth{\bE_{\theta\sim \pi}[\theta](e^{\lambda}-1) + \frac{(e^\lambda-1)^2 h^2}{8}} \\
&=\lambda(\bE_{\mu}\qth{X} - \bE_{\nu}\qth{X}) - (e^\lambda-\lambda-1)\bE_{\nu}\qth{X} - \frac{(e^\lambda-1)^2 h^2}{8} \\
&\ge \lambda(\bE_{\mu}\qth{X} - \bE_{\nu}\qth{X}) - (e^\lambda-\lambda-1)h - \frac{(e^\lambda-1)^2 h^2}{8}. 
\end{align*}
Let $\Delta:=\bE_{\mu}\qth{X} - \bE_{\nu}\qth{X}\in [-h,h]$. Since $0\le e^{\lambda} - \lambda - 1 \le \lambda^2$ holds for all $\lambda \in [-1,1]$, taking the supremum over $\lambda\in [-1,1]$ yields
\begin{align*}
\KL(\mu \| \nu) \ge \max_{\lambda\in [-1,1]} \bqth{\lambda \Delta - \lambda^2\bpth{h+\frac{h^2}{2}}} = \frac{\Delta^2}{4h+2h^2}, 
\end{align*}
where the maximum is attained at $\lambda^\star = \frac{\Delta}{2h+h^2}$, with $|\lambda^\star|\le \frac{|\Delta|}{2h}\le \frac{1}{2}$ thanks to $|\Delta|\le h$. 

\subsection{Regret lower bounds in the normal mean model}\label{append:regret_LB}
In this section, we establish the following lower bound on $\regS$ for case 3 of \Cref{thm:EB-gaussian}, thereby complete the discussions of known bounds above \Cref{cor:EB-gaussian}. 

\begin{thm}\label{thm:EB_regret_lb}
In a Gaussian location model with $\Theta=\sth{(\theta_1,\dots,\theta_n): \mu_p^w(G_n) \le C}$, where $\mu_p^w(G_n)$ is the weak $\ell_p$ norm of the empirical distribution $G_n$ defined in \eqref{eq:weak-ell-p-ball}, then
\begin{align*}
    \inf_{\widehat{\theta}}\regS(\widehat{\theta}) \ge c n^{\frac{1}{1+p}}\log^{-\frac{p}{2(1+p)}} n, 
\end{align*}
where $c>0$ is an absolute constant depending only on $p$ and $C$. 
\end{thm}
\begin{proof}
We apply Assouad's lemma to establish this lower bound. Define a sequence of real numbers $w_1,\dots,w_m$ as
\begin{align*}
    w_i = 4i \sqrt{\log n}, \quad i = 1, 2, \dots, m=\Big\lfloor c_0\frac{n^{1/(1+p)}}{\log^{\frac{p}{2(1+p)}}n} \Big\rfloor,
\end{align*}
with a small universal constant $c_0>0$. By simple algebra, this choice of $m$ ensures that the weak $\ell_p$ norm of the empirical distribution of $(w_1,\dots,w_m,0,\dots,0)$ is at most $C$. Next consider the following prior distribution $\pi$ on $(\theta_1,\dots,\theta_n)$: choose i.i.d. random bits $b_1,\dots,b_m\sim \Unif(\{0,1\})$, and set
\begin{align*}
\theta_i = w_i - b_i, \quad i=1,\dots,m, \qquad \theta_{m+1} = \dots = \theta_n = 0. 
\end{align*}
Since $\theta_i\le w_i$ for all $i\in [n]$, the condition on the weak $\ell_p$ norm is always satisfied. To lower bound the regret against the separable oracle, we need to 1) lower bound the MSE for any estimator $\widehat{\theta}$, and 2) upper bound the MSE for the best separable decision rule $\thetahatS$.

We start with the second target. Recall that $\thetahatS$ has the knowledge of $(\theta_1,\dots,\theta_n)$ but is restricted to take a separable form $\thetahatS_i = f(X_i)$ for a single function $f$. We construct the function $f$ as follows: 
\begin{align*}
    f(x) = \begin{cases}
        0 & \text{if } x\le 2\sqrt{\log n} \\
        \theta_i & \begin{aligned}
  &\text{if } (4i-2)\sqrt{\log n}<x\le(4i+2)\sqrt{\log n},\\
  &\text{and } 1\le i<m,
  \end{aligned}\\
        \theta_m & \text{if } x > (4m+2)\sqrt{\log n}
    \end{cases}. 
\end{align*}
We show that when $(\theta_1,\dots,\theta_n)\sim \pi$, for any fixed $i\in [n]$, the equality $\theta_i = f(X_i)$ holds with high probability. In fact, 
\begin{align*}
\bP(f(X_i)\neq \theta_i) \le \bP(|X_i - \theta_i|\ge 2\sqrt{\log n}-1) = O\pth{\frac{1}{n^2}}
\end{align*}
by Gaussian tail bounds. Consequently, 
\begin{align*}
\MSE(\theta,f) \le \sum_{i=1}^n \pth{\max_{i\in [n]} \theta_i}^2 \bP(f(X_i)\neq \theta_i) = O(1). 
\end{align*}
Finally, since $\thetahatS$ is the best separable decision rule (see \eqref{eq:separable_oracle}), we have
\begin{align}\label{eq:separable_upper}
\bE_{\theta\sim \pi}\qth{\MSE(\theta,\thetahatS)} \le \bE_{\theta\sim \pi}\qth{\MSE(\theta,f)} = O(1). 
\end{align}

Next we lower bound the Bayes MSE under the prior $\pi$. First, note that for $\theta^b$ and $\theta^{b'}$ determined by different bit strings $b,b'\in \{0,1\}^m$, we have the separation
\begin{align*}
\|\theta^b - \theta^{b'}\|^2 = \sum_{i=1}^m \indc{b_i \neq b_i'}. 
\end{align*}
Second, for neighboring $b$ and $b'$ (i.e. with $\sum_{i=1}^m \indc{b_i\neq b_i'}=1$), it is clear that
\begin{align*}
\KL(\calN(\theta^b,1) \| \calN(\theta^{b'},1)) = \frac{\|\theta^b - \theta^{b'}\|^2}{2} = \frac{1}{2} = O(1). 
\end{align*}
Therefore, by Assouad's lemma \cite[Theorem 2.12]{Tsy09}, 
\begin{align}\label{eq:bayes_lower}
\inf_{\widehat{\theta}}\bE_{\theta\sim \pi}[\MSE(\theta,\widehat{\theta})] = \Omega(m) = \Omega\pth{n^{\frac{1}{1+p}}\log^{-\frac{p}{2(1+p)}} n}. 
\end{align}

Finally, we conclude from \eqref{eq:separable_upper} and \eqref{eq:bayes_lower} that
\begin{align*}
\regS(\widehat{\theta}) &\ge \bE_{\theta\sim \pi}\qth{ \MSE(\theta,\widehat{\theta}) - \MSE(\theta,\thetahatS)} \\ &= \Omega\pth{n^{\frac{1}{1+p}}\log^{-\frac{p}{2(1+p)}} n},
\end{align*}
which is the claimed regret lower bound. 
\end{proof}